\documentclass{amsproc}

\usepackage{graphs}
  \usepackage{color}

\newtheorem{theorem}{Theorem}[section]
\newtheorem{lemma}[theorem]{Lemma}

\newtheorem{proposition}[theorem]{Proposition}

\theoremstyle{definition}
\newtheorem{definition}[theorem]{Definition}
\newtheorem{example}[theorem]{Example}

\theoremstyle{remark}
\newtheorem{remark}[theorem]{Remark}
\newtheorem{fact}[theorem]{Fact}

\numberwithin{equation}{section}

\def\rank{{\rm rank}}

\def\wt{\widetilde}
\def\ol{\overline}
\def\ov{\overline}
\def\tl{\widetilde}
\def\wh{\widehat}

\newcommand{\N}{{\mathbb N}}
\newcommand{\Z}{{\mathbb Z}}
\newcommand{\R}{{\mathbb R}}

\newcommand{\B}{{\mathcal B}}

\newcommand{\va}{\varphi}
\newcommand{\A}{\mathcal A}
\newcommand{\La}{\mathcal L}

\def\ov{\overline}

\newcommand{\om}{\omega}

\def\qed{\hfill$\square$}

\numberwithin{equation}{section}
\begin{document}

\title{Invariant measures for Cantor dynamical systems}

\author{S. Bezuglyi}
\address{Department of Mathematics, University of Iowa, Iowa City,
52242 IA, USA}
\email{bezuglyi@gmail.com}

\author{O. Karpel}
\address{Department of Mathematics, Institute for Low Temperature Physics, Kharkiv 61103, Ukraine}
\email{helen.karpel@gmail.com}

\curraddr{Faculty of Applied Mathematics, AGH University of Science and Technology, al. Mickiewicza 30, 30-059 Krakow, Poland}

\date{}

\dedicatory{Dedicated to the memory of our friend and colleague Sergiy Kolyada}


\begin{abstract}
This paper is a survey devoted to the study of probability 
and infinite ergodic invariant measures for aperiodic homeomorphisms of 
a Cantor set. We focus mostly  on the cases when a homeomorphism  has
 either a unique ergodic invariant measure or finitely many such measures 
 (finitely ergodic homeomorphisms). Since every Cantor dynamical system
 $ (X,T)$ can be realized as a Vershik map acting on the path space of a 
 Bratteli  diagram, we use combinatorial  methods developed in symbolic 
 dynamics and Bratteli  diagrams  during the last decade to study the 
 simplex of invariant measures. 
\end{abstract}

\maketitle
\tableofcontents

\section{Introduction}\label{introduction}

In this survey, we focus on an old classic problem of ergodic theory: for a
 given dynamical system $(X, \varphi)$, determine the set $M(X,\varphi)$ 
 of invariant  measures. In this generality, the problem is too complicated. 
 To make
 it more precise, we will consider only aperiodic Cantor dynamical systems,
 i.e., aperiodic homeomorphisms $\varphi$ of a Cantor set $X$. 
 There are many natural examples
 of such systems including subshifts in symbolic dynamics. 
We will discuss the significant progress which was made during the last 
decade in this direction. 

Because the problem of finding invariant measures for transformation arises 
in various areas of mathematics, we hope that this survey may be interesting 
not only for experts working in the ergodic theory but also for 
mathematicians who are interested in applications of these results. We 
included necessary definitions and formulated the most important facts to
make this text as much self-contained as possible. So that we begin with
the necessary background. 

Let $(X, \va )$ be a topological dynamical system, i.e., $\va$  is a 
homeomorphism of a compact metric space $X$. A Borel positive measure 
$\mu$ on $X$ is called \textit{invariant} if $\mu(\varphi(A)) = \mu(A)$ for
 any Borel set $A$. By the Kakutani-Markov theorem, such a measure 
 always exists. The set of all probability  invariant measures $M(X, \va)$ is 
 a Choquet simplex. Let $E(X, \va)$ 
denote the subset of extreme points of the simplex $M(X,\va)$. It is known 
that this set is formed by  ergodic measures for $\va$. By definition, a
 measure  $\mu$ is called \textit{ergodic} if $\varphi(A) = A$ implies that 
 either $\mu(A) =0$ or $\mu(X \setminus A) =0$. 
The cardinality of the set  $E(X, \va)$  can be either any positive integer, 
or $\aleph_0$, or continuum.  If $|E(X, \va)| =1$, then $\va$ is called 
\textit{uniquely ergodic}. If $|E(X, \va)| =k$, $k \in \N$, then $\va$ is called 
\textit{finitely ergodic}.
The reader can find more information in numerous books on ergodic theory 
and topological dynamics, we mention in this connection the books  \cite{CornfeldFominSinai1982},
\cite{Petersen1989}, and \cite{Walters1982}.

The question about a complete (or even partial) description of the simplex 
$M(X, \varphi)$  of invariant probability measures for $(X, \varphi)$ is one 
of the most important in ergodic theory. It has a long history and
many remarkable results.  The cardinality of the set of ergodic measures is 
an important invariant of dynamical
systems. The study of relations between the properties of the simplex 
$M(X, \varphi)$ and those of the dynamical system $(X, \varphi)$ is a hard
 and intriguing problem. There  is an extensive list of references regarding 
 this problem, we mention here only the books 
 \cite{Phelps2001}, \cite{Glasner2003}  
 and the papers \cite{Downarowicz1991},  \cite{Downarowicz2006},
 \cite{Downarowicz2008}  for further citations.
In particular, it is important to know conditions  under which a system
$(X, \varphi)$ is uniquely ergodic or has a finite number of ergodic measures.

We recall that the simplex $M(X, \varphi)$ plays an important role in 
the classification problems. In particular, it is a complete invariant for orbit
equivalence of minimal homeomorphisms of a Cantor set 
\cite{GiordanoPutnamSkau1995}.

The problem of finding invariant measures of a dynamical system 
$(X,\varphi)$ looks rather 
vague in general setting.  There are very few universal results that can be
 applied to  an arbitrary   homeomorphism $\varphi$.  
 A very productive idea is to replace $(X, \varphi)$ by an isomorphic  model 
 $(X_B, \varphi_B)$ for which the computation of invariant measures is
 more transparent.  
To study invariant measures for a Cantor system $(X, \varphi)$, we will work 
with  {\em Bratteli diagrams}, the object that is widely used for constructions  
of  transformation models in various dynamics, see Section 
\ref{Section_basics_CS_BD} for definitions. 
It is difficult to overestimate the significance of Bratteli diagrams for the 
theory of dynamical systems. A class of graduated infinite graphs, later called 
Bratteli diagrams, was originally introduced by Bratteli \cite{Bratteli1972} 
in his breakthrough article on the classification of approximately finite 
 $C^*$-algebras.

It turned out that the ideas developed by Vershik in the ergodic theory
\cite{Vershik1981}, \cite{Vershik1982} found their application in Cantor 
dynamics. It was proved in \cite{HermanPutnamSkau1992}, that any minimal
homeomorphisms $\varphi$ of a Cantor set $X$  can be represented as a homeomorphism $\varphi_B$ 
(called \textit{Vershik map}) acting on the path space $X_B$ of a Bratteli diagram $B$. 
 The dynamical systems obtained in this way are called  
 \textit{Bratteli-Vershik dynamical systems.} Later on, this approach was realized for 
 non-minimal Cantor dynamical systems \cite{Medynets2006} and Borel
  automorphisms of  a standard Borel space \cite{BezuglyiDooleyKwiatkowski2006}. 

The literature devoted to Cantor dynamical systems is very extensive. 
We do not plan to discuss many interesting directions such as the 
classification of homeomorphisms up to orbit equivalence, dimension groups, 
the interplay of Cantor dynamical systems and  $C^*$-algebras, etc.
The  reader, who is interested in this subject, can  be referred to the 
recent surveys and books
\cite{Skau2000}, \cite{Putnam2010}, \cite{Durand2010}, 
\cite{BezuglyiKarpel2016}, \cite{Putnam2018} and the research papers 
\cite{GiordanoPutnamSkau1995},  \cite{GiordanoPutnamSkau1999}, \cite{GiordanoMatuiPutnamSkau2010}, 
\cite{DurandHostSkau1999} (more references can be found in the cited 
surveys). 

The main reason why Bratteli diagrams are convenient to use for the study 
of homeomorphisms $\varphi:  X\to X$  is the fact that various properties of
 $\varphi$ become more transparent when one deals with the corresponding 
Bratteli-Vershik dynamical systems. This observation is related first of all to 
$\varphi$-invariant measures and 
their supports, to minimal components of $\varphi$, structure of 
$\varphi$-orbit, etc.
In particular, the study of an ergodic $\varphi$-invariant measure $\mu$ 
 is  reduced, roughly speaking, to the computation of the values of $\mu$ on
cylinder subsets in the path space of the corresponding Bratteli diagram. 
In other words, the \textit{structure} of a Bratteli diagram determines
completely the invariant measures. In this case we should speak about 
the invariance with respect to \textit{the tail equivalence relation} because 
there are Bratteli diagrams that do not admit Vershik maps.  We 
emphasize the difference between simple and non-simple Bratteli 
diagrams in this context.   For an aperiodic homeomorphism $\varphi$,
 the simplex 
$M(X,\varphi)$ may contain the so called ``regular'' infinite measures, i.e., 
the  infinite $\sigma$-finite measures that take finite (nonzero) values on 
some  clopen sets.  

We give one more important observation about Bratteli diagrams. They
can be used to construct homeomorphisms of a Cantor set with 
\textit{prescribed} properties. For instance, it is easy to build a diagram
that has exactly $k$ ergodic invariant measures. 

A similar picture occurs in symbolic dynamics. Let $(X, S)$ be a subshift, i.e.,
$X$ is a shift invariant closed subset of $\A^{\Z}$ where $\A$ is a finite 
alphabet. Then combinatorial structure of sequences from the set $X$ can
be used to determine invariant measures. Boshernitzan's results  
about  the number of ergodic measures for  a minimal subshift give  
the bounds in terms of the complexity function (see 
\cite{Boshernitzan1984}. We see a clear similarity between the application of
Bratteli diagrams and complexity functions to estimate the number of
ergodic measures. This is a reason why we included the recent results 
extending Boshernitzan's approach (see \cite{CyrKra2019}, 
\cite{DamronFickenscher2017}, \cite{DamronFickenscher2019} and the 
references therein).
 
In the current paper, we focus on the following problem: how  determine
the number of ergodic measures for a given Cantor dynamical system.   
We distinguish three classes of dynamical systems: uniquely ergodic, 
finitely ergodic, and ``infinitely ergodic'' systems. This problem was  
considered in symbolic dynamics for minimal subshifts by many authors
(see the references in Sections \ref{Section_basics_CS_BD}, 
\ref{Section_Uniq_Erg}, and \ref{finrank}).

The outline of the paper is as follows. In Section \ref{Section_basics_CS_BD},
we give necessary definitions and facts that are used below in the main
text. The key concepts are Bratteli diagrams (ordered, simple, non-simple,
stationary, finite rank, etc), subshifts, complexity functions. Section 
\ref{Section_Inv_measures} contains a description  of the simplex 
of invariant measures in terms of incidence matrices. We also  discuss 
the problem of measure extension from a subdiagram. In other words,
the proved results clarify conditions that would guarantee finiteness of measures invariant with respect to the tail equivalence relation. 
In Section \ref{Section_Uniq_Erg}, we collected results about uniquely 
ergodic dynamical systems. These results are formulated either in terms
of complexity functions or in terms of Bratteli diagrams. We understand that 
the variety of uniquely ergodic transformations is very vast, and the
included results have been chosen to illustrate the discussed methods.
In the next section, we consider the results about dynamical systems that
have finitely many ergodic measures. For us, the most important sources of 
examples are stationary and finite rank Bratteli diagrams. In the last section,
Setion \ref{inf_rank_section}, we consider a class of Bratteli diagrams
that have countably many ergodic invariant measures. In the paper, the reader 
will find a big 
number of explicit examples to picturize the principal theorems.

\section{Basics on Cantor dynamics and Bratteli diagrams}\label{Section_basics_CS_BD}

This section contains the basic definitions and facts about topological
(in particular, Cantor) dynamical systems. Most of the definitions can be 
found in the well known books on topological and symbolic dynamics, we
refer to \cite{LindMarkus1995}, \cite{Kitchens1998}, \cite{Kurka2003}. 

\subsection{Cantor dynamical systems}

A{\em Cantor set (space)} $X$ is  a zero-dimensional compact metric space
 without isolated points. The topology on $X$ is generated by a countable 
 family of clopen subsets. All such Cantor sets are homeomorphic. 

For a homeomorphism $T : X \to X$, denote $Orb_T(x) := \{T^n (x) \; |\; n
 \in \mathbb{Z}\}$; the set $Orb_T(x)$ is called the {\em orbit} of $x\in X$ 
under action of $T$ (or simply $T$-orbit). We consider here mostly {\em 
aperiodic} homeomorphisms $T$, i.e., for every $x$ the set $Orb_T(x)$  is 
countably infinite.

A homeomorphism $T : X \to X$ is called {\em minimal} if for every $x \in X$
 the orbit $Orb_T(x)$ is dense. Any (aperiodic) homeomorphism $T$ of a
Cantor  set has a {\em minimal component}: this is a $T$-invariant 
closed non-empty subset $Y$ of $X$ such that $T|_Y$ is minimal on $Y$.

There are several natural notions of equivalence for  Cantor dynamical 
systems. We 
give the definitions of conjugacy and orbit equivalence for single 
homeomorphisms of Cantor sets.

\begin{definition}
Let $(X,T)$ and $(Y, S)$ be two  Cantor systems. Then

(1) $(X,T)$ and $(Y,S)$ are \textit{conjugate} (or \textit{isomorphic}) if 
there exists a homeomorphism $h \colon X \rightarrow Y$ such that $h \circ 
T = S \circ h$.

(2) $(X,T)$ and $(Y,S)$ are \textit{orbit equivalent} if there exists a 
homeomorphism $h \colon X \rightarrow Y$ such that $h(Orb_T(x)) = 
Orb_S(h(x))$ for every $x \in X$. In other words, there exist functions $n, m 
\colon X \rightarrow \mathbb{Z}$ such that for all $x \in X$, $h \circ T(x) = 
S^{n(x)} \circ h(x)$ and $h \circ T^{m(x)} = S \circ h(x)$. 

\end{definition}

Let $\B$ be the Borel sigma-algebra generated by clopen subsets of $X$.
We consider only Borel positive measures on $(X, \B)$. A measure $\mu$
is called probability (finite) if $\mu(X) =1$ ($\mu(X) < \infty$).  
Similarly, $\mu$ is infinite if $\mu(X) =\infty$. In the latter, we assume that
$\mu$ is a sigma-finite measure.  
Given a Cantor dynamical system $(X,T)$, a Borel measure $\mu$ on $X$ is 
called $T$-\textit{invariant} if $\mu(TA) = \mu(A)$ for any $A\in \B$. 
A measure $\mu$ is called \textit{ergodic} with respect to $T$ if, for any 
$T$-invariant  set $A$, either $\mu(A) = 0$ or $\mu(X \setminus A) =0$.

Let $M(X,T)$ be the set of all $T$-invariant probability  measures. It is well
 known that  
$M(X,T)$ is a Choquet simplex whose extreme points are exactly  
$T$-ergodic measures.  Denote by $E(X,T)$ the set of extreme points 
(ergodic measures) in $M(X,T)$. If $M(X,T) =\{\mu\}$, then $T$ is called 
{\em uniquely ergodic}.  Clearly, that in this case  $|E(X,T)| =1$ where $|
\cdot |$ denotes the cardinality of a set.  If $|E(X,T)| =k  < \infty$, we say
 that $(X,T)$ is \textit{finitely ergodic.} 

It can be easily seen that if two systems,  $(X, T)$ and $(Y,S)$, are 
conjugate (or orbit equivalent), then $M(X,T) = M(Y,S)$.

\subsection{Languages on finite alphabets and complexity}
\label{subsect languages complexity}
In this subsection, we recall the definitions from symbolic dynamics. 
This material can be found in many books, see e.g. \cite{LindMarkus1995}. 

We first recall several definitions from symbolic dynamics. Let $\mathcal A$
be a finite \textit{alphabet}, then a \textit{word} $w = a_1\cdots a_k$ in 
this alphabet is a concatenation of letters  $a_i$ in $\mathcal A$. The 
length $|w|$ is the number of letters in $w$. Let $\mathcal A^n$ denote
the set of words over $\mathcal A$ of length $n$. Then,
by $\mathcal A^* = \bigcup_{n=1}^\infty \A^n$, we denote the set of all 
finite nonempty words. It is said that a word $w = a_1\cdots a_k$ occurs in 
a word  $u = b_1 \cdots b_s$ if $a_1 = b_m,  ... , a_k = b_{m+k -1}$. The
 word $w$ is called a \textit{subword} (or \textit{factor}) of $u$.

A language $\mathcal L$ can be determined in the 
abstract setting as follows.

\begin{definition}\label{def language} 
A set $\mathcal L$ of finite words on  an
 alphabet $\mathcal A$ is called a \textit{language} if:
 
(i) $\mathcal A \subset \mathcal L$, 

(ii) for any word $w$ from  
$\mathcal L$, all subwords $w' $ of $ w$ belong to $\mathcal L$ 
(the language is \textit{factorial)}; 

 (iii) for any word $w \in \mathcal L$, there exist letters $a$ and $b$ such
 that $awb \in\mathcal L$ (the language is \textit{extendable)}. 
 \end{definition}
 
Let $\mathcal L_n = \mathcal A^n \cap \mathcal L$ denote the set  words 
 in the language $\mathcal L$ of length $n$.
 
 A language $\La$ is called \textit{recurrent} if for any $u, v \in \La$
 there exists a word $w \in \La$ such that $uwv \in \La$, and 
 $\La$ is called \textit{uniformly} \textit{recurrent} if for every $u \in \La$
there exists $m \in \N$ such that $u$ is a subword in every $w\in \La_m$.
We consider  \textit{aperiodic} languages only ($\La$ is periodic if for every 
word $w= a_1 \cdots a_{|w|}$ there 
exists $p \in \N$ such that $a_i = a_{i + p}$  where $1 \leq |w| - p$.

The notion of a language is naturally arisen in symbolic dynamical systems.
We first note that for every infinite sequence $\om\in \A^{\N}$ of symbols 
from $\A$, one can define the language $\La(\om)$ determined by 
$\omega$ as the family of all finite subwords that occur in $\om$. 

More generally, one  can define the \textit{language of a subshift} 
$(X, S)$ where  $S : \A^{\Z} \to \A^{\Z}$ denote the (left) shift, and $X $
 is a closed  $S$-subset of $\A^{\Z}$.
 For a subshift $(X, S)$, we define the language $\La(X)$ of 
$(X, S)$ as the set of all finite words that occur in the sequences $x$ 
from $X$. Clearly, $\La(X)$ is a factorial and extendable language.
 Conversely, if a language $\La$ on a finite alphabet $\A$ is defined, then 
 there exists
a subshift $(X_{\La}, S)$ whose language coincides with $\La$. Indeed, the 
set $X_{\La}$ is now determined by those sequences from $\A^{\Z}$
whose finite subwords belong to $\La$. It is obvious that 
$$
\La(X_{\La}) = \La.
$$  
In other words, the map $(X, S) \to   \La(X)$ is a bijection from the set 
of subshifts to the set of non-empty factorial and extendable languages.

The dynamical properties of subshifts $(X, S)$ can be represented in terms 
of the corresponding languages. For example, the dynamical system 
$(X, S)$ is \textit{minimal} if and
only if the language $\La(X)$ is uniformly recurrent. In this case, the 
language $\La(X)$ coincides with $\La(\om)$ where $\omega$ is an 
arbitrary point (sequence) from $X$. If one fixes a sequence $\om$ in 
$\A^{\N}$, then the language $\La(\om)$ defines a subshift denoted
by $(X_\om, S)$.
\\

For every language $\La$, we define the \textit{symbolic complexity
 function}  $p_{\La}(n) : \N \to \N$ by setting
$$
p_{\La}(n ) = |\La_n|
$$
where $|\cdot|$ stands for the cardinality of a set. If the language $\La$ is 
defined by a sequence $u \in \A^{\Z}$, then the corresponding complexity
function is denoted by $p_u(n)$. Clearly, the complexity function is 
non-decreasing.

Let $(X, S)$ be a minimal subshift on a finite alphabet $\A$. Then
the complexity function $p_X(n)$ can be defined either as that of the
corresponding language $\La(X)$ or that of an infinite sequence $\om \in 
X$. In both cases, these functions are the same.

The complexity functions have been
studied extensively in many papers devoted to languages and symbolic
dynamical systems, see, e.g. the survey  \cite{Ferenczi1999} and the 
bibliography therein.  We mention 
here several results about the complexity functions of dynamical systems.

\begin{fact}
(i) \cite{Ferenczi1996} Let  $(X_u, S)$ and $(X_v, S)$ be symbolic 
dynamical systems defined by uniformly recurrent sequences $u$ and $v$ 
from $\A^{\N}$. If $(X_u, S)$ and $(X_v, S)$ are topologically conjugate. 
Then there exists an integer $c$ such that, for all $n > c$,
$$
p_u(n -c ) \leq p_v(n) \leq p_u(n +c).
$$

(ii) \cite{CovenHedlund1973} 
Let $\omega$ be a sequence in $\A^{\N}$. Then $\omega$ is 
ultimately periodic if and only if there exists $n \geq 1$ such that 
$p_\omega(n) \leq n$ if and only if there exists $n \geq 1$ such that 
$p_\omega(n) = p_\omega(n+1)$.

(iii) \cite{Pansiot1984} Let $\zeta : \A \to \A^*$ be a primitive substitution.
Then  the complexity function $p_u(n)$ of the sequence $u = \zeta(u)$ is 
sublinear, i.e.,  there exists $C$, a positive constant, such that
$p_u(n) \leq Cn$, for $n \geq 1$.

Moreover, the set of differences $p_u(n+1) - p_u(n)$ is bounded 
\cite{Cassaigne1996}.

(iv) $p_u(m + n) \leq p_u(m) p_u(n) $, and the limit $\lim_n n^{-1}
\log p_u(n)$ is the topological entropy of a sequence. 

(v) There are sequences $u$ such that $p_u(n)= n+1$; they are  called   
Sturmian sequences.

\end{fact}

\subsection{Ordered Bratteli diagrams and Vershik maps}
 A {\it Bratteli diagram} is an infinite graph $B=(V,E)$ such that the vertex
set $V =\bigcup_{i\geq 0}V_i$ and the edge set $E=\bigcup_{i\geq 0}E_i$
are partitioned into disjoint subsets $V_i$ and $E_i$ where

(i) $V_0=\{v_0\}$ is a single point;

(ii) $V_i$ and $E_i$ are finite sets, $\forall i \geq 0$;

(iii) there exist $r : E \to V$ (range map $r$) and $s : E \to V$ (source map 
$s$) such that $r(E_i)= V_{i+1}$, $s(E_i)= V_{i}$, and
$s^{-1}(v)\neq\emptyset$, $r^{-1}(v')\neq\emptyset$ for all $v\in V$
and $v'\in V\setminus V_0$.

The set of vertices $V_i$ is called the $i$-th level of the
diagram $B$. A finite or infinite sequence of edges $(e_i : e_i\in E_i)$
such that $r(e_{i})=s(e_{i+1})$ is called a {\it finite} or {\it infinite
path}, respectively. For $m<n$, $v\, \in V_{m}$ and $w\,\in V_{n}$, let 
$E(v,w)$ denote the set of all paths $\overline{e} = (e_{1},\ldots, e_{p})$
with $s(\ol e) = s(e_{1})=v$ and $r(\ol e) = r(e_{p})=w$. 
For a Bratteli diagram $B$,
let $X_B$ be the set of infinite paths starting at the top vertex $v_0$.
 We endow $X_B$ with the topology generated by cylinder sets
$[\ol e]$ where $\ol e= (e_0, ... , e_n)$,  $n \in \mathbb N$, and
$[\ol e]:=\{x\in X_B : x_i=e_i,\; i = 0, \ldots, n\}$.
 With this topology, $X_B$ is a 0-dimensional compact metric space.
By assumption, we will consider only such Bratteli diagrams $B$ for which 
$X_B$ is a {\em Cantor set}, that is $X_B$ has no isolated points.  

Given a Bratteli diagram $B$, the $n$-th {\em incidence matrix}
$F_{n}=(f^{(n)}_{v,w}),\ n\geq 0,$ is a $|V_{n+1}|\times |V_n|$
matrix such that $f^{(n)}_{v,w} = |\{e\in E_{n+1} : r(e) = v, s(e) = w\}|$
for  $v\in V_{n+1}$ and $w\in V_{n}$.  Every vertex $v \in V$ is connected 
with $v_0$ by a finite path, and the set of $E(v_0, v) $ of all such paths is 
finite. If $h^{(n)}_v = |E(v_0, v)|$, then, for all $n \geq 1$, we have
\begin{equation}\label{formula for heights}
h_v^{(n+1)}=\sum_{w\in V_{n}}f_{v,w}^{(n)}h^{(n)}_w \ \mbox{or} \ \
h^{(n+1)}=F_{n}h^{(n)}
\end{equation}
where $h^{(n)}=(h_w^{(n)})_{w\in V_n}$. The numbers $h_w^{(n)}$ are 
usually called \textit{heights} (the terminology comes from using the 
Kakutani-Rokhlin partitions to build a Bratteli-Vershik system for a 
homeomorphism of a Cantor space, see below).

We define the following important classes of Bratteli diagrams:

\begin{definition}\label{rank_d_definition} Let $B$ be a Bratteli diagram.
\begin{enumerate}
\item
We say that $B$ has \textit{finite rank} if  for some $k$, $|V_n| \leq k$ for 
all $n\geq 1$.
\item
 We say that a finite rank diagram $B$ has \textit{rank $d$} if  $d$ is the 
 smallest integer such that $|V_n|=d$ infinitely often.
\item
We say that $B$ is {\em simple} if for any level
$n$ there is $m>n$ such that $E(v,w) \neq \emptyset$ for all $v\in
V_n$ and $w\in V_m$. Otherwise, $B$ is called {\em non-simple}.
\item
 We say that $B$ is  \textit{stationary} if $F_n = F_1$  for all $n\geq 2$.
\end{enumerate}
\end{definition}

Let $x =(x_n)$ and $y =(y_n)$  be two paths from $X_B$. It is said that 
$x$ and $y$ are {\em tail equivalent} (in symbols,  $(x,y) \in \mathcal R$)  if 
there exists some $n$ such that $x_i = y_i$ for all $i\geq n$. Since $X_B$ 
has no isolated points, the $\mathcal R$-orbit of any point $x\in X_B$ is 
infinitely countable. The diagrams with infinite $\mathcal R$-orbits are called 
{\em aperiodic}. Note that a Bratteli diagram is simple if the tail equivalence 
relation $\mathcal R$ is minimal.

In order to illustrate the above definitions, we give an example of a 
nonsimple Bratteli diagram.

\unitlength = 0.5cm
 \begin{center}
 \begin{graph}(17,15)
\graphnodesize{0.4}
 \roundnode{V0}(8.5,14)
\roundnode{V11}(0.0,9) \roundnode{V12}(4.25,9)
\roundnode{V13}(8.5,9) \roundnode{V14}(12.75,9)
\roundnode{V15}(17.0,9)
\roundnode{V21}(0.0,5) \roundnode{V22}(4.25,5)
\roundnode{V23}(8.5,5) \roundnode{V24}(12.75,5)
\roundnode{V25}(17.0,5)
\roundnode{V31}(0.0,1) \roundnode{V32}(4.25,1)
\roundnode{V33}(8.5,1) \roundnode{V34}(12.75,1)
\roundnode{V35}(17.0,1)
\edge{V11}{V0} \edge{V12}{V0} \edge{V13}{V0} \edge{V14}{V0}
\edge{V15}{V0}
\edge{V21}{V11}\bow{V21}{V11}{0.06} \edge{V21}{V12} \edge{V22}{V11} \edge{V22}{V12}
\edge{V23}{V13} \edge{V23}{V14} \edge{V24}{V13} \edge{V24}{V14}
\edge{V25}{V11}\bow{V25}{V11}{0.06} \edge{V25}{V13} \bow{V25}{V15}{0.06}
\bow{V25}{V15}{-0.06} \edge{V25}{V15}
\edge{V31}{V21} \edge{V31}{V22} \bow{V31}{V22}{0.06}\edge{V32}{V21} \edge{V32}{V22}
\edge{V33}{V23}\bow{V33}{V23}{0.06} \edge{V33}{V24} \edge{V34}{V23} \edge{V34}{V24}
\edge{V35}{V21} \edge{V35}{V23} \bow{V35}{V25}{0.06}
\bow{V35}{V25}{-0.06} \edge{V35}{V25}

\end{graph}
\centerline{.\ .\ .\ .\ .\ .\ .\ .\ .\ .\ .\ .\ .\ .\ .\ .\ }
\\
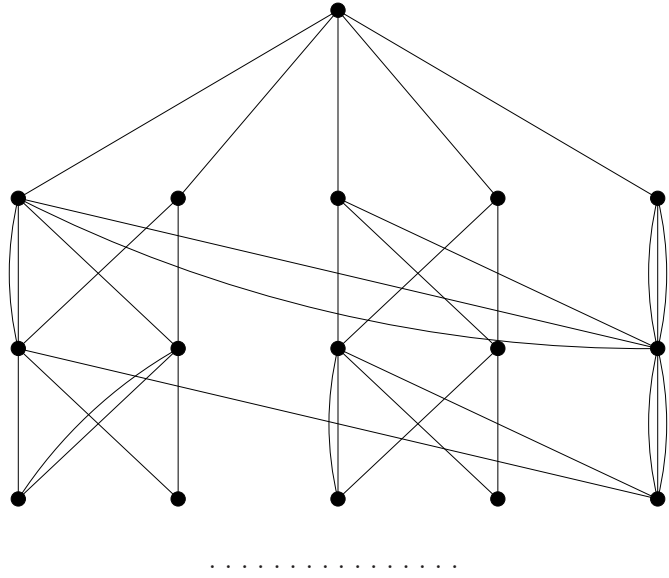
Fig. 1. Example of a Bratteli diagram

 \end{center}
\medskip

This diagram is a non-simple  finite rank Bratteli diagram that has exactly two 
minimal components (they are clearly seen).

We will  constantly use the \textit{telescoping} procedure for a Bratteli 
diagram:

\begin{definition} \label{telescoping_definition}
Let $B$ be a Bratteli diagram, and $n_0 = 0 <n_1<n_2 < \ldots$ be a 
strictly increasing sequence of integers. The {\em telescoping of $B$ to $
(n_k)$} is the Bratteli diagram $B'$, whose $k$-level vertex set $V_k'$ is 
$V_{n_k}$ and whose incidence matrices $(F_k')$ are defined by
\[F_k'= F_{n_{k+1}-1} \circ \ldots \circ F_{n_k},\]
where $(F_n)$ are the incidence matrices for $B$.
\end{definition}

Roughly speaking, in order to telescope a Bratteli diagram, one takes a 
subsequence of levels $\{n_k\}$ and considers the set $E(n_k, n_{k+1})$  of 
all finite paths between the  levels $\{n_k\}$ and $\{n_{k+1}\}$ as edges of 
the new diagram.
In particular, a Bratteli diagram $B$ has rank $d$ if and only if there is
a telescoping $B'$ of $B$ such that $B'$ has exactly $d$ vertices
at each level. When telescoping diagrams, we often do not specify to which 
levels $(n_k)$ we telescope, because it suffices to know that such a 
sequence of levels exists.

To avoid consideration of some trivial cases, we will assume that the following 
\textit{convention} always holds: {\em our Bratteli diagrams are not unions of 
two or more disjoint subdiagrams.}

The concept of an ordered Bratteli diagram is crucial for the existence of 
dynamics on the path space of a Bratteli diagram.

\begin{definition}\label{order_definition} A Bratteli diagram $B=(V,E) $ is 
called {\it ordered} if a linear order `$>$' is defined on every set  $r^{-1}(v)
$, $v\in
\bigcup_{n\ge 1} V_n$. We denote by $\om$  the corresponding partial
order on $E$ and write $(B,\om)$ when we consider $B$ with the ordering 
$\om$. Let $\mathcal O_{B}$ denote the set of all orders on $B$.
\end{definition}

Every $\omega \in \mathcal O_{B}$ defines the  \textit{lexicographic}
order on the set $E(k,l)$ of finite paths between
vertices of levels $V_k$ and $V_l$:  $(e_{k+1},...,e_l) > (f_{k+1},...,f_l)$
if and only if there is $i$ with $k+1\le i\le l$, such that $e_j=f_j$ for 
$i<j\le l$ and $e_i> f_i$.
It follows that, given $\om \in \mathcal O_{B}$, any two paths from 
$E(v_0, v)$
are comparable with respect to the lexicographic order generated by 
$\om$.  If two infinite paths are tail equivalent, and agree from the vertex 
$v$ onwards, then we can compare them by comparing their initial segments 
in $E(v_0,v)$. Thus, $\om$ defines a partial order on $X_B$, where two 
infinite paths are comparable if and only if they are tail equivalent.

   \begin{definition}
Let $(B, \omega)$ be an ordered Bratteli diagram.  We call a finite or infinite
 path $e=(e_i)$ \textit{ maximal (minimal)} if every $e_i$ is maximal
  (minimal) amongst the edges from the set $r^{-1}(r(e_i))$.
\end{definition}

Denote by $X_{\max}(\om)$ and $X_{\min}(\om)$ the sets of all maximal and
minimal infinite paths in $X_B$, respectively. It is not hard to see that
$X_{\max}(\om)$ and $X_{\min}(\om)$ are \textit{non-empty closed subsets} of
$X_B$; in general, $X_{\max}(\om)$ and $X_{\min}(\om)$ may have
interior points. For a finite rank Bratteli diagram $B$, the sets $X_{\max}(\om)$
and $X_{\min}(\om)$ are always finite for any $\om$, and if $B$ has rank $d$,
then each of them have at most $d$ elements
(\cite{BezuglyiKwiatkowskiMedynets2009}). For an aperiodic Bratteli diagram $B$, we see that   $X_{\max}(\om)\cap X_{\min}(\om) =\emptyset$.

We say that an  ordered Bratteli diagram $(B, \om)$ is  \textit{properly ordered} if the sets $X_{\max}(\om)$ and $X_{\min}(\om)$ are singletons. A Bratteli diagram is called \textit{regular} if the set
of maximal paths and set of minimal paths have empty interior.

\begin{definition}\label{VershikMap}
Let $(B, \omega)$ be an ordered Bratteli
diagram. We say that $\varphi = \varphi_\omega : X_B\rightarrow X_B$
is a {\it (continuous) Vershik map} if it satisfies the following conditions:

(i) $\varphi$ is a homeomorphism of the Cantor set $X_B$;

(ii) $\varphi(X_{\max}(\omega))=X_{\min}(\omega)$;

(iii) if an infinite path $x=(x_0,x_1,\ldots)$ is not in $X_{\max}(\omega)$,
then $\varphi(x_0,x_1,\ldots)=(x_0^0,\ldots,x_{k-1}^0,\overline
{x_k},x_{k+1},x_{k+2},\ldots)$, where $k=\min\{n\geq 1 : x_n\mbox{
is not maximal}\}$, $\overline{x_k}$ is the successor of $x_k$ in
$r^{-1}(r(x_k))$, and $(x_0^0,\ldots,x_{k-1}^0)$ is the minimal path
in $E(v_0,s(\overline{x_k}))$.
\end{definition}

If $\om$ is an ordering on $B$, then one can always define the map
$\varphi_0$ that maps $X_B \setminus X_{\max}(\om)$ onto $X_B
\setminus X_{\min}(\om)$ according to (iii) of Definition
\ref{VershikMap}. The question about the existence of the Vershik map is
equivalent to that of an extension of  $\varphi_0 : X_B \setminus
X_{\max}(\om) \to X_B \setminus X_{\min}(\om)$ to a homeomorphism
of the entire set $X_B$.  If $\om$ is a proper ordering, then
$\varphi_\om$ is a homeomorphism. In particular any simple Bratteli diagram 
admits a Vershik map. 

\begin{definition}\label{Good_and_bad}  Let $B$ be a Bratteli diagram
$B$.  We say that an ordering $\om\in \mathcal O_{B}$ is \textit{perfect}
if $\om$ admits a Vershik map $\varphi_{\om}$ on $X_B$.
Denote by  $\mathcal P_B $ the set of
all perfect orderings  on $B$.
\end{definition}

We observe that for a regular Bratteli diagram with an order $\om$, the
 Vershik map $\varphi_\om$, if it exists, is
defined in a unique way. Also, a necessary condition for $\om\in
\mathcal P_{B}$ is that $|X_{\max}(\om)|=|X_{\min}(\om)|$. Given 
$(B,\om)$ with $\om \in \mathcal P_B$, the uniquely defined  system 
$(X_B, \varphi_\om)$ is called  a {\em Bratteli-Vershik} or {\em adic} 
system.

We can summarize the above definitions and results in the following  
statement.

\begin{theorem}  Let $B= (V,E, \om)$ be an ordered Bratteli diagram with 
a perfect order
$\om  \in \mathcal P_B$. Then there exists an aperiodic homeomorphism 
(Vershik map) $\varphi_{\om}$ acting on the path space $X_B$ according to 
Definition \ref{VershikMap}. The homeomorphism $\varphi_\om$ is minimal if 
and only if $B$ is simple.
\end{theorem}

The pair $(X_B, \varphi_\om)$ is called the {\em Bratteli-Vershik dynamical 
system}.

The simplest example of a Bratteli diagram is an {\em odometer}. Any 
odometer can be realized as a Bratteli diagram $B$ with $|V_n| =1$ for all 
$n$. Then any order on $B$ is proper and defines the Vershik map.

It is worth noticing that a general Bratteli diagram may have a rather 
complicated structure. In particular, the tail equivalence relation may have 
uncountably many minimal components or, in other words, uncountably many 
simple subdiagrams that do not have connecting edges.

The ideas developed in the papers by Vershik \cite{Vershik1981}, 
\cite{Vershik1982}, where sequences of refining measurable partitions of a 
measure space were used to construct a realization of an ergodic 
automorphisms of a measure space, turned out to be very fruitful for finding a 
model of any minimal homeomorphism $T$ of a Cantor set $X$. In 
\cite{HermanPutnamSkau1992}, Herman, Putnam, and Skau found an explicit 
construction that allows one to define an ordered simple Bratteli diagram $B = 
(V, E, \om)$ such that $T$ is conjugate to the corresponding Vershik map $
\varphi_\omega$. The authors used the existence of the first return time 
map to any clopen set to build the nested sequence of Kakutani-Rohklin 
partitions and the corresponding ordered Bratteli diagram. Since this 
construction is described in many papers (not only in 
\cite{HermanPutnamSkau1992}), we will not  give the details here referring 
to the original paper and \cite{Durand2010} for detailed explanation. The 
case of aperiodic Cantor system is much subtler and was considered in 
\cite{BezuglyiDooleyMedynets2005} and \cite{Medynets2006}. 

Let $(X, T)$ be an aperiodic Cantor system. A closed subset $Y$ of $X$ is 
called {\em basic} if (1) $Y \cap T^iY = \emptyset, i\neq 0$, and (2) every 
clopen neighborhood $A$ of $Y$ is a complete $T$-section, i.e., $A$ meets 
every $T$-orbit at least once. This means that every point from $A$ is 
recurrent. It is clear that if $T$ is minimal then every point of $X$ is a basic 
set.  It was proved in \cite{Medynets2006} that every aperiodic Cantor 
system $(X, T)$ has a basic set. This is a crucial step in the proof of 
the following theorem.

\begin{theorem} \cite{Medynets2006}\label{existenceBD}
Let $(X, T)$ be a Cantor aperiodic system with a basic set $Y$. There exists 
an ordered Bratteli diagram  $B = (V,E,\om)$ such that $(X, T)$ is conjugate 
to a Bratteli-Vershik dynamical system $(X_B, \varphi_\om)$. The 
homeomorphism implementing the conjugacy between $T$ and $\varphi_\om
$ maps the basic set $Y$ onto the set $X_{\min}(\om)$ of all minimal paths 
of $X_B$. The equivalence class of $B$ does not depend on a choice of 
$\{\xi(n)\}$ with the property $\bigcap_n B(\xi(n)) =Y$.

\end{theorem}

Is the converse theorem true? In case of a \textit{simple Bratteli diagram, the 
answer is  obviously affirmative:} there exists a proper order $\om$ on any 
simple Bratteli diagram $B$ so that $(X_B, \varphi_\om)$ is a minimal Cantor 
system. \textit{For general non-simple Bratteli diagrams the answer is
 negative}. The first example of a Bratteli diagram that does not admit 
 a Vershik map was found in \cite{Medynets2006}. A systematic study of this
 problem is given in 
\cite{BezuglyiKwiatkowskiYassawi2014}, 
\cite{BezuglyiYassawi2017}, \cite{JanssenQuasYassawi2017}, see also
\cite{BezuglyiKarpel2016}. 

Perfect orderings were also studied in \cite{DownarowiczKarpel2018, DownarowiczKarpel2019}. The approach there was a bit different: the starting point was the abstract compact invertible zero-dimensional system $(X,T)$ and the aim was to find an ordered regular Bratteli diagram $B = (V,E,\om)$ with the perfect ordering $\om$ such that $(X_B, \varphi_\om)$ is topologically conjugate to $(X,T)$. Regular perfectly ordered Bratteli diagrams are called \textit{decisive}. The following theorem holds:
\begin{theorem}[\cite{DownarowiczKarpel2019}]\label{main}
A (compact, invertible) zero-dimensional system $(X,T)$ is topologically conjugate to a decisive Bratteli-Vershik system $(X_B, \varphi_\om)$ if and only if the set of aperiodic points of $(X,T)$ is dense, or its closure misses one periodic orbit.
\end{theorem}
The proof uses Krieger's Marker Lemma~\cite{Boyle1983} and representation of $(X,T)$ as an array system. Also in~\cite{Shimomura2018} a non-trivial Bratteli-Vershik model is build for every compact metric zero-dimensional dynamical system.

\section{Invariant measures on Bratteli diagrams}\label{Section_Inv_measures}


Since any aperiodic Cantor dynamical system $(X,T)$ admits a realization as a Bratteli-Vershik dynamical system (see Section~\ref{Section_basics_CS_BD}), the study of $T$-invariant measures is reduced to the case of measures defined on the path space of a Bratteli diagram. The advantage of this approach is based on the facts that (i) any such a measure is completely determined by its values on cylinder sets of $X_B$, and (ii) there are simple and explicit formulas for measures of cylinder sets. Especially transparent this method works for stationary and finite rank Bratteli diagrams, simple and non-simple ones \cite{BezuglyiKwiatkowskiMedynetsSolomyak2010}, \cite{BezuglyiKwiatkowskiMedynetsSolomyak2013}.

It is worth pointing out that the study of measures on a Bratteli diagram is a more general problem than that in  Cantor dynamics. This observation follows from the existence of Bratteli diagrams that do not support any continuous dynamics on their path spaces which is compatible with the tail equivalence relation. The first example of such a Bratteli diagram was given  in \cite{Medynets2006}; a more comprehensive coverage of this subject can be found in \cite{BezuglyiKwiatkowskiYassawi2014} and \cite{BezuglyiYassawi2017} (see also \cite{JanssenQuasYassawi2017} and \cite{BezuglyiKarpel2016}). If a Bratteli diagram does not admit a Bratteli-Vershik homeomorphism, then we have to work with the {\em tail equivalence relation} $\mathcal R$ on $X_B$ and study measures invariant with respect to $\mathcal R$.

\subsection{Simplices, stochastic incidence matrices, examples}\label{Subsection_simplices}
In this subsection, we show that the set of all probability invariant measures
on a Bratteli diagram corresponds to the inverse limit of a decreasing  sequence of convex sets.
Let $\mu$ be a Borel probability non-atomic $\mathcal R$-invariant
  measure on $X_B$ (for brevity, we will use the term
  ``measure on $B$'' below). We denote the set of all such measures
   by $\mathcal{M}_1(B)$ and by $\mathcal{E}_1(B)$ the set of all probability ergodic invariant measures.
   The   fact that $\mu$ is an
   $\mathcal R$-invariant measure
 means that  $\mu([e]) = \mu([e'])$ for any two cylinder sets
 $e,e' \in  E(v_0, w)$, where $w \in  V_n$ is an arbitrary vertex, and
 $n \geq1$.   Since any measure on $X_B$ is completely determined
  by its  values  on  clopen (even cylinder) sets, we conclude that in order to
  define an  $\mathcal{R}$-invariant measure $\mu$, one needs to
   know the  sequence   of vectors  $\ov p^{(n)} = (p_w^{(n)} :
   w \in V_n), n  \geq 1$,  such that $p_w^{(n)} =  \mu([e(v_0,w)])$
    where  $e(v_0,w)$ is a finite path
 from $E(v_0,w)$. It is clear that, for $w \in V_n$,
\begin{equation}\label{eq path extension}
[e(v_0, w)] = \bigcup_{e(w, v), v \in V_{n+1}} [e(v_0, w), e(w, v)]
\end{equation}
so that $[e(v_0, v)] \subset [e(v_0, w)]$.
Relation (\ref{eq path extension}) implies that
\begin{equation}\label{formula for measures}
\tl F_n^T \ov p^{(n+1)} = \ov p^{(n)}, \ n\geq 1,
\end{equation}
where  $\tl F_n^T$ denotes the transpose of the matrix $\tl F_n$.
The entries of the vectors $\ov p^{(n)}$ can
 be also found by the formula
$$
p_w^{(n)} = \frac{\mu(X_w^{(n)})}{h_w^{(n)}},
$$
where
\begin{equation}\label{eq tower X}
X_w^{(n)} = \bigcup_{e \in E(v_0,w)} [e], \ \ w \in V_n.
\end{equation}
The clopen set $X_w^{(n)}$ is called a \textit{tower}, since it is the tower in the Kakutani-Rokhlin
  partition that corresponds to the vertex $w$
   (see Section~\ref{Section_basics_CS_BD}).  The
  measure of this tower is
 \begin{equation}\label{eq tower measure}
  \mu(X_w^{(n)}) = h_w^{(n)}p_w^{(n)} =: q_w^{(n)}.
 \end{equation}
 Denote $\ov q^{(n)} = (q_w^{(n)} :    w \in V_n), n  \geq 1$.

Because $\mu(X_B) =1$, we see that,  for any $n >1$,
$$
\sum_{w \in V_n} h_w^{(n)}p_w^{(n)} = \sum_{w \in V_n} q_w^{(n)} = 1.
$$
We can obtain a formula similar to (\ref{formula for measures}), but for $\ov q^{(n)}$ instead of $\ov p^{(n)}$ and stochastic incidence matrices $F_n$ instead of usual incidence matrices $\tl F_n$. 
The entries of the row stochastic incidence matrix $F_n$ are defined by the formula
\begin{equation}\label{stoch_inc_matr}
f_{vw}^{(n)} = \frac{\tl f_{vw}^{(n)} h_w^{(n)}}{h_v^{(n+1)}}.
\end{equation}

\begin{example}[Equal row sums (ERS) Bratteli diagrams]
\label{ExampleERS1} In this example, we compute the stochastic incidence matrices for a class
of Bratteli diagrams that have the so called \textit{equal row sum (ERS}) property.
A Bratteli diagram $B$ has the ERS property if there exists a sequence of natural numbers $(r_n)$ such that the incidence matrices $(\wt F_n)$ of $B$ satisfy the condition
  $$
  \sum_{w \in V_n} \tl f_{v,w}^{(n)} = r_n
  $$
  for every $v \in V_{n+1}$. It is known that Bratteli-Vershik
   systems with the ERS property can serve as models for
   Toeplitz subshifts   (see \cite{GjerdeJohansen2000}).
In particular, we have $\wt F_0 = \ov h^{(1)} = (r_0, \ldots, r_0)^T$.
 It follows from (\ref{formula for heights}) that, for ERS Bratteli diagrams,
  $h^{(n)}_w = r_0 \cdots r_{n-1}$ for every $w \in V_n$.
Hence we have for all $n \geq 1$, $w \in V_n$ and $v \in V_{n+1}$:
$$
f_{vw}^{(n)} = \frac{\tl f_{vw}^{(n)}}{r_n}.
$$
\end{example}

In general, it is difficult to compute the elements of the matrix $F_n$ explicitly because the terms $h_w^{(n)}$ used in the formula (\ref{stoch_inc_matr}) are the entries of the product of matrices. In Section~\ref{Section_Uniq_Erg}, the reader can find Examples~\ref{ex 3.2} and \ref{non-simple-stat-ex} of Bratteli diagrams, for which stochastic incidence matrices are computed explicitly. Some of the  results about the exact number of ergodic invariant measures for a diagram are formulated in terms of $\ov q^{(n)}$ and $F_n$ (see Sections~\ref{finrank} and~\ref{inf_rank_section}). It is easy to prove 
the following lemma.
\begin{lemma}\label{formula_for_qn}
Let $\mu$ be a probability measure on the path
 space $X_B$ of a Bratteli diagram $B$. Let $(F_n)$ be a sequence
 of corresponding stochastic incidence matrices. Then, for every $n
  \geq 1$, the vector $\ov q^{(n)} = ( q^{(n)}_v : v \in
   V_n)$ (see (\ref{eq tower measure})) is a  probability vector  such that
 \begin{equation}\label{eq stoch matrix and measures}
 F_n^T \ov q^{(n+1)} = \ov q^{(n)}, \ \ \ \ n  \geq 1.
 \end{equation}
\end{lemma}

We see that the formula in (\ref{eq stoch matrix and measures}) is a
 necessary
condition for a sequence of vectors $(\ov q^{(n)})$ to be defined by an
  invariant  probability measure. It turns out that the  converse
   statement is  true, in general. We formulate below Theorem
    \ref{Theorem_measures_general_case}, where all
 $\mathcal{R}$-invariant measures are explicitly
described.

Using the definition of stochastic incidence matrix $F_n$
 (see~(\ref{stoch_inc_matr})) and Lemma~\ref{formula_for_qn},  we
 define a decreasing sequence of convex polytopes $\Delta_m^{(n)}$,
 $n,m \geq 1$, and the limiting convex sets $\Delta_{\infty}^{(n)}$. They
 are used  to describe the set  $\mathcal{M}_1(B)$ of all probability
 $\mathcal{R}$-invariant measures on $B$. Namely, denote
$$
\Delta^{(n)} := \{(z^{(n)}_w)_{w \in V_n}^T : \sum_{w \in V_n}
z^{(n)}_w = 1\ \mbox{ and } \ z^{(n)}_w   \geq 0, \  w \in V_n \}.
$$
The sets $\Delta^{(n)}$ are standard simplices in the space
$\mathbb{R}^{|V_n|}$ with $|V_n|$
 vertices $\{\ov e^{(n)}(w) : w \in V_n\}$, where
$\ov e^{(n)}(w) = (0, ... , 0, 1, 0,...0)^T$ is the standard basis vector, i.e.
$e^{(n)}_u(w) = 1$ if and only if $u = w$.
Since $F_n$ is a stochastic matrix, we have the obvious property
 $$
 F_{n}^T(\Delta^{(n+1)}) \subset \Delta^{(n)}, \quad n \in \mathbb{N}.
 $$

Let $\mu$ be  a probability $\mathcal{R}$-invariant measure $\mu$
on $X_B$ with values $q_w^{(n)}$ on the towers $X_w^{(n)}$.
Then $(q_w^{(n)} : w \in V_n)^T$ lies in  the simplex
 $\Delta^{(n)}$.
Set
\begin{equation}\label{eq def Delta^n_m}
\Delta_m^{(n)} = F_n^T \cdot \ldots \cdot
F_{n+m-1}^T(\Delta^{(n+m)})
\end{equation}
for $m = 1,2,\ldots$
Then we see that
$$
\Delta^{(n)} \supset \Delta^{(n)}_1 \supset \Delta^{(n)}_2 \supset \ldots
$$
Denote
\begin{equation}\label{eq convex set Delta}
\Delta_{\infty}^{(n)} = \bigcap_{m=1}^{\infty}\Delta_m^{(n)}.
\end{equation}
It follows from  (\ref{eq def Delta^n_m}) and  (\ref{eq convex set Delta})
 that
\begin{equation}\label{eq F_n relates Delta}
F_n^T(\Delta_{\infty}^{(n +1)}) = \Delta_{\infty}^{(n)}, \quad n \geq 1.
\end{equation}

The next theorem, that was proved in
\cite{BezuglyiKwiatkowskiMedynetsSolomyak2010}, describes
all $\mathcal{R}$-invariant probability measures. 

\begin{theorem}[\cite{BezuglyiKwiatkowskiMedynetsSolomyak2010}, \cite{BezuglyiKarpelKwiatkowski2018}]
\label{BKMS_measures=invlimits}
 \label{Theorem_measures_general_case}  Let $B = (V,E)$ be a
 Bratteli diagram with the sequence of stochastic incidence matrices
 $(F_n)$, and let $\mathcal M_1(B)$ denote the set of
 $\mathcal R$-invariant probability measures on the path space
 $X_B$.

(1) If  $\mu \in \mathcal M_1(B)$, then the probability vector
$$
\ov q^{(n)}= (\mu(X_w^{(n)}))_{w\in V_n},
$$
where $X_w^{(n)}$ is defined in (\ref{eq tower X}), satisfies the following conditions for $n\geq 1$:

(i) $\ov q^{(n)} \in  \Delta_{\infty}^{(n)}$,

(ii) $F^{T}_n\ov q^{(n+1)} = \ov q^{(n)}$.

Conversely, suppose that $\{\ov q^{(n)}\}$ is a sequence of
non-negative probability vectors such that, for every
 $\ov q^{(n)}= (q^{(n)}_w)_{w\in V_n} \in  \Delta^{(n)}$ ($n\geq 1$),
 the condition  (ii) holds. Then the vectors  $\ov q^{(n)}$ belong to
 $\Delta_{\infty}^{(n)}$, $n \in \mathbb N$, and   there  exists a uniquely
  determined $\mathcal R$-invariant  probability measure $\mu$ such that
   $\mu(X_w^{(n)})= q_w^{(n)}$   for $w\in V_n, n \in \mathbb N$.

(2) Let $\Omega$ be the subset of the infinite product
$\prod_{n\geq 1} \Delta^{(n)}_\infty$ consisting of sequences
$(\ov q^{(n)})$ such that $F^{T}_n\ov q^{(n+1)} = \ov q^{(n)}$. Then
 the map
$$
\Phi : \mu \ \mapsto \ (\ov q^{(n)}) : \mathcal M_1(B)  \mapsto
\Omega,
$$
is an affine isomorphism. Moreover, $\Phi(\mu)$ is an extreme point of
$\Omega$ if and only if $\mu$ is ergodic.

(3) Let $B$ be a Bratteli diagram of rank $K$. Then the number of ergodic invariant
 measures on $B$ is
 bounded  above by $K$ and bounded below by the dimension of  the finite-dimensional simplex
 $\Delta^{(1)}_\infty$.
\end{theorem}

\begin{remark}
(a) From Theorem \ref{Theorem_measures_general_case} it follows that the set
$\mathcal{M}_1(B)$ can be identified with the inverse limit of the
sequence $(F_n^T, \Delta^{(n)}_{\infty})$.
In general, the set $\Delta_\infty^{(n)}$ is a convex subset of the
$(|V_n| - 1)$-dimensional simplex $\Delta^{(n)}$. In some cases,
 which will be considered in Section~\ref{finrank}, the set
 $\Delta_\infty^{(n)}$ is a finite-dimensional simplex itself.

(b)  In Theorem~\ref{Theorem_measures_general_case} part (2), the set
$\mathcal{M}_1(B)$ can be affinely isomorphic to the set
$\Delta_{\infty}^{(1)}$. For instance, it happens when all stochastic
incidence matrices are square non-singular matrices of the same dimension
$K \times K$ for some $K \in
\mathbb{N}$. This case will be considered in Section~\ref{finrank}.

(c) The procedure of telescoping (see Definition~\ref{telescoping_definition})
preserves the set of invariant measures; hence we can apply it when
necessary without loss of generality.
\end{remark}

In order to find all ergodic invariant measures on a Bratteli diagram, we will study the
number of  extreme points of $\Delta^{(n)}_{\infty}$ for every $n$.

%

Let
\begin{equation}\label{matrix_G_def}
G_{(n+m,n)} = F_{n+m} \cdots F_{n}
\end{equation}
for $m \geq 0$ and $n \geq 1$.
Denote the elements of $G_{(n+m,n)}$ by $(g_{uw}^{(n+m,n)})$, where
$u \in V_{n+m+1}$ and $w\in V_{n}$.
The sets $\Delta_m^{(n)}, m \geq 0$, defined in
(\ref{eq def Delta^n_m}),
form a decreasing sequence of convex polytopes in $\Delta^{(n)}$.
The vertices of $\Delta_m^{(n)}$
are some (or all) vectors from the set $\{\ov g^{(n+m,n)}(v) :
v \in V_{n+m+1}\}$, where we denote
$$
\ov g^{(n+m,n)}(v) = (g_{w}^{(n+m,n)}(v))_{w \in V_n} =
G_{(m+n,n)}^T(\ov e^{(n+m+1)}(v)).
$$
Obviously, we have the relation
\begin{equation}\label{2.3}
\ov g^{(n+m,n)}(v) = \sum_{w \in V_n}g_{vw}^{(n+m,n)}\ov e^{(n)}(w).
\end{equation}
Let $\{\ov y^{(n,m)}(v)\}$ be the set of all vertices of
$\Delta_m^{(n)}$.
Then $\ov y^{(n,m)}(v) = \ov g^{(n+m,n)}(v)$ for $v$ belonging
 to some subset 
 of $V_{n+m+1}$.

We observe the following fact. Every vector $\ov q^{(n)} $ from
 the set $\Delta_{\infty}^{(n)}$ can be  written in the standard basis as
$$
\ov q^{(n)} = \sum_{w\in V_n}q^{(n)}_w \ov e^{(n)}(w).
$$
 It turns out that the numbers $q_v^{(n+m+1)}$, $v \in V_{n+m+1}$,
  are the coefficients in the convex decomposition of $\ov q^{(n)}$
  with respect to vectors $\ov g^{(n+m,n)}(v)$.

\begin{proposition}[\cite{BezuglyiKarpelKwiatkowski2018}]\label{decomposition_q_n_by_q_n+m}
Let $\mu \in \mathcal{M}_1(B)$, and $q^{(n)}_w = \mu(X^{(n)}_w)$
($w \in V_n$) for all $n \in \mathbb{N}$. Then
\begin{equation}\label{2.3ab}
\ov q^{(n)} = \sum_{v \in V_{n+m+1}} q_v^{(n+m+1)} \ov
g^{(n+m,n)}(v).
\end{equation}
In particular,
$$
\ov q^{(1)} = \sum_{v \in V_{m+1}} q_v^{(m+1)} \ov g^{(m+1,1)}(v).
$$
\end{proposition}

\begin{remark}
For every $n \geq 1$,  define
$$
\Delta^{(n), \varepsilon}_{\infty} := \bigcup_{\ov q \in
 \Delta_{\infty}^{(n)}} B(\ov q, \varepsilon),
$$
where $B(\ov q, \varepsilon)$ is the ball of radius $\varepsilon > 0$
centered at  $\ov q \in \mathbb{R}^{|V_n|}$. Here the metric is
defined by the Eucleadian norm $||\cdot ||$ on $\mathbb R^{|V_n|}$.
Fix any natural numbers  $n$ and $m$.   Let $\Delta_{m}^{(n)}$ be defined
as above. It can be proved straightforwardly that if $\ov q^{(n,m)} \in \Delta_{m}^{(n)}$ for infinitely many $m$
and $\ov q^{(n,m)} \rightarrow \ov q^{(n)}$ as $m \to \infty$,   then
$\ov q^{(n)} \in
\Delta_{\infty}^{(n)}$. Moreover, for every $\varepsilon > 0$ there exists $m_0 = m_0(n,\varepsilon)$
 such that $\Delta_{m}^{(n)} \subset \Delta^{(n),
 \varepsilon}_{\infty}$ for all $m \geq m_0$.
\end{remark}

The next statement shows that vertices of the limiting convex set
$\Delta^{(n)}_{\infty}$ can be obtained as limits of sequences of vertices
 of convex polytopes. 

\begin{lemma}[\cite{BezuglyiKarpelKwiatkowski2018}]\label{remar3}
Fix $n \in \mathbb{N}$ and $\varepsilon >0$.  Let
$\Delta^{(n)}_{\infty}$ and
$\Delta_{n+m+1}^{(n)}$, 
be defined as above for
any $m \in \mathbb{N}$.  Then, for every vertex $\ov y \in
\Delta^{(n)}_{\infty}$ there exists $m_0 = m_0(n,\varepsilon)$ such that,
 for all $m \geq m_0$, one can find a vertex $\ov y^{(n,m)}(v) \in
  \Delta_{n+m+1}^{(n)}$ 
satisfying the property
$$
 ||\ov y - \ov y^{(n,m)}(v)|| < \varepsilon.
 $$
\end{lemma}


\subsection{Subdiagrams and measure extension (finite and infinite measures)}
In this subsection, we study supports of ergodic invariant measures on arbitrary Bratteli diagrams in terms of subdiagrams. By a Bratteli subdiagram, we mean a Bratteli diagram $\ov B$ that can be obtained from $B$ by removing some vertices and edges from each level of $B$. Then $X_{\ov B}\subset X_B$. We will consider two extreme cases of Bratteli subdiagrams: vertex subdiagram (when we fix a subset of vertices at each level and take all existing edges between them) and edge subdiagram (some edges are removed from the initial Bratteli diagram but the set of vertices is preserved). It is clear that an arbitrary subdiagram can be obtained as a combination of these cases.

Take a subdiagram $\ol B$ and consider the set $X_{\ol B}$ of all infinite paths whose edges belong to $\ol B$. As a rule, objects related to a subdiagram $\ol B$ are denoted by barred symbols. Let $\widehat X_{\ol B} := \mathcal R(X_{\ol B})$ be the subset of paths in $X_B$ that are tail equivalent to paths from $X_{\ol B}$.
Let $\ov \mu$ be a probability measure on $X_{\ol B}$ invariant with respect to the tail equivalence relation defined on $\ol B$. Then $\ov \mu$ can be canonically extended to the measure $\widehat {\ov \mu}$ on the space $\widehat X_{\ol B}$ by invariance with respect to $\mathcal R$ \cite{BezuglyiKwiatkowskiMedynetsSolomyak2013, AdamskaBezuglyiKarpelKwiatkowski2016}. If we want to extend $\widehat{\ov \mu}$ to the whole space $X_{B}$, we set $\widehat {\ov \mu} (X_B \setminus \widehat{X}_{\ol B}) = 0$.

This subsection is devoted to answering the following questions:

(A) Given a subdiagram $\ov B$ of $B$ and an ergodic measure $\mu$ on $X_B$, under what conditions on $\ov B$ the subset $X_{\ov B}$ has positive measure $\mu$ in $X_B$?

(B) Let $\nu$ be a measure supported by the path space $X_{\ov B}$ of a  subdiagram $\ov B \subset B$. Then $\nu$ is extended to the subset $\mathcal R(X_{\ov B})$ by invariance with respect to the tail equivalence relation $\mathcal R$. Under what conditions is $\nu(\mathcal R(X_{\ov B}))$ finite (or infinite)?

In this subsection, we keep the following notation:
$\overline{X}_v^{(n)}$ stands for the tower in a subdiagram $\overline{B}$ that is determined by  a vertex $v$ of $\ol B$. Thus,  we consider the paths in $\overline{X}_v^{(n)}$ that contain edges from $\overline{B}$ only. Let $\overline{h}_v^{(n)}$ be the  height of the tower $\ol X_v^{(n)}$. 
The following theorem gives criteria for finiteness of the measure extension. 

\begin{theorem}[\cite{BezuglyiKarpelKwiatkowski2015}]\label{thm from BKK}
Let $B$ be a Bratteli diagram with the sequence of incidence matrices $\{\tl F_n\}_{n=0}^\infty$ and corresponding stochastic matrices $\{F_n\}_{n=0}^\infty$. Let $\ov B$ be a vertex subdiagram of $B$ defined by the sequence of subsets  $\{W_n\}_{n = 0}^\infty, \ W_n \subset V_n$. Suppose that $\ov \mu$ is  a probability  $\ol{\mathcal R}$-invariant measure on $X_{\ov B}$. Then the following properties are equivalent:
\begin{eqnarray*}
\widehat{\overline{\mu}}(\wh X_{\ov B}) < \infty & \Longleftrightarrow & \sum_{n=1}^\infty \sum_{v \in W_{n+1}} \sum_{w \notin W_n} \tl f_{v,w}^{(n)} h_w^{(n)} \ov p_v^{(n+1)} < \infty\\
& \Longleftrightarrow & \sum_{n=1}^{\infty} \sum_{v \in W_{n+1}} \wh{\ov\mu}(X_{v}^{(n+1)})\sum_{w \notin W_n} f_{v,w}^{(n)}
< \infty\\
& \Longleftrightarrow & \sum_{i=1}^{\infty} \left(\sum_{w\in W_{i+1}} h_w^{(i+1)} \ov p_w^{(i+1)} - \sum_{w\in W_{i}} h_w^{(i)} \ov p_w^{(i)}\right)< \infty.
\end{eqnarray*}
\end{theorem}

The analogue of Theorem~\ref{thm from BKK} can be proved also for edge subdiagrams (see   \cite{AdamskaBezuglyiKarpelKwiatkowski2016}). The following proposition gives a sufficient condition of the finiteness of the measure extension (more necessary and sufficient conditions can be found in~\cite{BezuglyiKarpelKwiatkowski2015, AdamskaBezuglyiKarpelKwiatkowski2016}).

\begin{proposition}[\cite{BezuglyiKarpelKwiatkowski2015}]\label{1.6b}
Let $B$ be a Bratteli diagram with the sequence of stochastic incidence matrices $\{F_n\}_{n=0}^\infty$, and let $\ov B$ be its subdiagram defined by a
 sequence of vertices $W_n$.
If
$$
\sum_{n=1}^\infty \max_{v \in W_{n+1}} \left(\sum_{w \notin W_n} f_{vw}^{(n)}\right) < \infty,
$$
then any tail invariant probability measure $\ov \mu$ on $X_{\ov B}$ extends to
a finite invariant measure $\wh {\ov \mu}$ on $\wh X_{\ov B}$.
\end{proposition}


The following theorem gives a necessary and sufficient condition for a  subdiagram $\ol B$ of $B$ to have a path space of zero measure in $X_B$.   Though the theorem is formulated for a vertex subdiagram, the statement remains  true also for any edge subdiagram $\ol B$.

\begin{theorem}[\cite{AdamskaBezuglyiKarpelKwiatkowski2016}]\label{main}
Let $B$ be a simple Bratteli diagram, and let $\mu$ be any probability ergodic measure on $X_B$. Suppose that  $\ov B$ is a vertex subdiagram of $B$ defined by a sequence  $(W_n)$ of subsets of $V_n$. Then $\mu(X_{\ov B}) = 0$ if and only if for all $\varepsilon > 0$ there exists $n = n(\varepsilon)$ such that for all $w \in W_n$ one has
\begin{equation}\label{thin set condition}
\frac{\ov h_w^{(n)}}{h_w^{(n)}} < \varepsilon.
\end{equation}
\end{theorem}

In fact, Theorem \ref{main} states that if  a subdiagram $\ol B$ satisfies  (\ref{thin set condition}), then $X_{\ol B}$ has measure zero with respect to every ergodic invariant measure, that is the set $X_{\ol B}$ is {\em thin} according to the definition from \cite{GiordanoPuntamSkau2004}. The following result is a corollary of Theorem \ref{main}:

\begin{theorem}[\cite{AdamskaBezuglyiKarpelKwiatkowski2016}]\label{corol_infty}
Let $\ol B$ be a subdiagram of $B$ such that $X_{\ol B}$ is a thin subset of $X_B$.  Then for any probability invariant measure $\ov \mu$ on $\ov B$ we have $\widehat{\overline{\mu}}(\wh X_{\ov B}) = \infty$.
\end{theorem}

\begin{remark}
There are many papers, where invariant measures for various Bratteli diagrams are studied. For instance, in~\cite{FrickPetersen2008, PetersenVarchenko2010} the authors consider ergodic invariant probability measures on a Bratteli diagram of a special form, called an Euler graph; the combinatorial properties of the Euler graph are connected to those of Eulerian numbers. The authors of~\cite{FrickOrmes2013} study spaces of invariant measures for a class of dynamical systems which is called polynomial odometers. These are adic maps on regularly structured Bratteli diagrams and include the Pascal and Stirling adic maps as examples. K. Petersen~\cite{Petersen2012} considers ergodic invariant measures on a Bratteli-Vershik dynamical system, which is based on a diagram whose path counts from the root are the Delannoy numbers.

We would like to mention also the interesting paper by Fisher \cite{Fisher2009} where various properties of Bratteli diagrams and measures are discussed.
\end{remark}

\section{Uniquely ergodic Cantor dynamical systems}
 \label{Section_Uniq_Erg}

In this section, we consider a number of results
 about uniquely ergodic Cantor dynamical systems. We are not trying to 
 mention all existing classes of uniquely ergodic homeomorphisms. In the case 
 of symbolic systems, we discuss the results related mostly to the complexity
function (Subsection 4.1). More general approach using Bratteli diagrams is 
considered in Subsection 4.2.

\subsection{Minimal uniquely ergodic homeomorphisms in symbolic 
dynamics}

In this subsection, we partially use some statements formulated and proved 
in \cite{FerencziMonteil2010} and in \cite{Boshernitzan1984},  
\cite{Boshernitzan1985}, and 
\cite{Boshernitzan1992}.

\begin{definition} Let $\om = (\om_i)$ be an infinite sequence in $\A^{\N}$.
It is said that the infinite sequence $\om$  has \textit{uniform frequencies} 
if, for every factor $w$ of $\om$ and any $a \in \A$, 
$$
\frac{|\om_k\cdots \om_{k+n}|_a}{n+1} \to f_w(\om)\quad (n \to \infty)
$$
uniformly in $k$. 
(Here $|u|_a$  denotes the number of occurrences  of the letter $a$ in 
the word $u$)
\end{definition}

The following result follows immediately from the ergodic theorem.

\begin{fact}[Folklore]
(i)  Let  $(X, S)$ be a subshift, and let $\mu$ be an $S$-invariant ergodic 
 measure. Then, for $\mu$-a.e. $\om \in X$ and for any finite word $w$ in
$\La(X)$, the frequency $f_w(\om)$ exists and is equal to $\mu([w])$ 
where $[w]$ denotes the corresponding cylinder subset of $X$.

(ii) A subshift  $(X_\om,  S)$ is uniquely ergodic if, and
only if, the sequence $\om$  has uniform frequencies.
\end{fact}

As one of our goals is to discuss relations  between the complexity functions
and the number of ergodic measures, we recall the Boshernitzan's results
 about uniquely ergodic subshifts. In fact, Boshernitzan proved several
  impressive results on the cardinality of the set of ergodic measures for  minimal subshifts which are 
 based on a careful study of the growth of the complexity functions.
The case of 
finite ergodicity is considered below in Section \ref{finrank}.

\begin{theorem}[\cite{Boshernitzan1984}]
Let $p_X(n)$ denote the complexity function of a minimal subshift $(X, S)$
over a finite alphabet. If either 
$$
\limsup_{n\to \infty} \frac{p_X(n)}{n} < 3,
$$
or 
$$
\liminf_{n\to \infty} \frac{p_X(n)}{n} = \alpha < 2,
$$
then $(X, S)$ is uniquely ergodic. 
\end{theorem}

Let $X$ be a compact metric space, $\B$ the Borel sigma-algebra, 
and $\mu$ a Borel probability measure on $\B$. Suppose  $T : X  \to X$
is a measurable map preserving the measure $\mu$.  A point $x \in X$ is
called  a \textit{generic point} for the measure $\mu$ if for every 
continuous function $f : X \to \R$,
$$
\lim_{n \to \infty} \frac{1}{n} \sum_{i=0}^{n-1} f(T^ix) = \int_X f \;d\mu.
$$
The measure $\mu$  is \textit{generic} if it has a generic point. 

By the pointwise ergodic theorem, if $\mu$ is ergodic, then almost every 
point is generic. 
However, there are generic measures which are not ergodic. An example of
such a measure is given in \cite{ChaikaMasur2015} where an interval
 exchange transformation has a generic non-ergodic measure. 

\begin{theorem} \cite{Glasner2003} Let $(X,T)$ be a minimal system. 
Then $(X, T)$ is 
uniquely ergodic if and only if every point $x$ in $X$ is generic for some
 measure in $M(X,T)$.
\end{theorem}

\begin{remark}
In this remark we point out several classes of uniquely ergodic Cantor
 dynamical systems.

 (1) \textit{Substitution dynamical systems}. Let $\A$ be a finite
 alphabet, and let $\tau : \A \to \A^*$ be a substitution. The language of
 the substitution consists of all words which are seen in 
 $\bigcup_{n=0}^\infty  \tau^n(\A)$. The corresponding subshift 
 $(X_\tau, S)$ is called a \textit{substitution dynamical system}. The
 literature on these dynamical systems
is very extensive, we refer to \cite{Queffelec2010},  \cite{Fogg2002}. A substitution 
$\tau$ is called \textit{primitive} if for any $a, b \in \A$ there exists
$n\in \N$ such that $|\tau^n(b)|_a \geq 1$. The corresponding 
substitution dynamical system is uniquely ergodic.

In Section \ref{finrank} below,  we consider \textit{aperiodic} 
(non-minimal) substitution systems.
 The situation 
with invariant measures is different. They may have finitely many egodic 
invariant probability measures and finitely many infinite ergodic invariant 
measures as well.

(2) \textit{Linearly recurrent dynamical systems}. 
Let $\A$ be an alphabet, and let $\om$ be a sequence from  $\A^{\N}$ 
with the language $\La(\om)$. For a word $u \in \La(\om)$, we call 
a word $w$ a return word to $u$ in $\om$ if $wu$ belongs to $\La(\om)$,
$u$ is a prefix of $wu$, and $u$ has exactly two occurrences in $wu$ (we 
follow \cite{Durand1998}, more general approach to the notion of return 
words  is given in \cite{Durand2010}).
Denote by $\mathcal R_{\om, u}$ the set of return
words to $u$ of $\om$. 
When $\om$ is a uniformly recurrent sequence from $\A^{\N}$, then  
all $u \in \La(\om)$ the set $\mathcal R_{\om, u}$ is finite.

It is said  that a sequence $\om$ is\textit{ linearly recurrent} 
 (with constant $K \in \N$) if it is uniformly recurrent and if for all $u \in 
 \La(\om)$ and all $w \in \mathcal R_{\om, u}$  we
have $|w| \leq  K|u|$. A subshift $(X, S)$ is  called \textit{linearly recurrent}
 (with constant $K$) if it is minimal and contains a linearly recurrent
  sequence (with constant $K$). In fact, for any $x, y \in X$, we have
$\mathcal R_{x, u} =\mathcal R_{y, u}$.

As proved in \cite{Durand2000} (see also \cite{Durand2010}), 
linearly recurrent subshifts are uniquely
 ergodic. This result can be deduced  from \cite{Boshernitzan1992} or 
 proved directly.
One more important fact that relates linearly recurrent subshifts and Bratteli 
diagrams is proved in \cite{Durand2010}). It states that such subshifts have
an expansive Bratteli-Vershik representation whose incidence matrices belong
to a finite set.

\end{remark}

\subsection{Finite rank Bratteli diagrams and general case}
\label{subsect_finrank_ue}

In this subsection, we discuss the results on unique ergodicity of Bratteli
diagrams. It is worth recalling  that these results describe Cantor dynamical 
systems which are represented by the corresponding Bratteli diagrams. We 
give a criterion and sufficient conditions for the unique ergodicity of a Bratteli 
diagram $B$ of arbitrary rank, in other words, we discuss the case when the 
space $\mathcal M_1(B)$ is a singleton.

We first begin with a class of Bratteli diagrams that  have an exact finite 
rank. 
\begin{definition}\label{def exact rank}
It is said that a finite rank Bratteli diagram  has an {\it exact finite rank} if 
there is a finite invariant measure $\mu$ and a constant $\delta >0$ such 
that after a telescoping $\mu(X_w^{(n)})\geq \delta$ for all levels $n$ and 
vertices $w$.
\end{definition}

The following result shows that the Vershik map on the path space of an 
exact finite rank diagram cannot be strongly mixing independently of the 
ordering.

\begin{theorem}[\cite{BezuglyiKwiatkowskiMedynetsSolomyak2013}]
\label{TheoremAbsenceStrongMixing} Let $B = (V,E,\omega)$
be an ordered simple Bratteli diagram of exact finite rank and $\mu$ is as
in Definition \ref{def exact rank}. 
\\
(1) The diagram $B$ is uniquely ergodic and $\mu$ is the unique invariant 
measure.
\\
(2) Let $\varphi_{\omega}:X_B\rightarrow X_B$ be the Vershik map 
defined by the
order $\omega$ on $B$ ($\varphi_{\omega}$ is not necessarily continuous 
everywhere).
 Then the dynamical system $(X_B,\mu,\varphi_{\omega})$ is not strongly 
 mixing with respect to the unique invariant measure $\mu$.
\end{theorem}

On the other hand, it is proved in the same paper that for the so-called ``left-
to-right'' ordering, the Vershik map is not strongly mixing on all finite rank 
diagrams.

\begin{remark} Theorem \ref{TheoremAbsenceStrongMixing} can be viewed
 as an analogue of 
a result from \cite{Boshernitzan1992}. In more details, let $(X, S)$ be a 
minimal subshift on a finite alphabet, and let $\mu$ be a probability 
$S$-invariant measure. Set 
$$
\varepsilon (n) = \min\{ \mu([w])  : w \in \La_n(X)\}
$$
where $[w]$ is the cylinder subset of $X$ defined by the word $w$. If 
$$
\lim_{n \to\infty} n \varepsilon(n) = 0,
$$
then the subshift $(X, S)$ is not uniquely ergodic. 
\end{remark}

In what follows, we focus on the following problem: find conditions on 
the (stochastic) incidence matrices under which the diagram is uniquely 
ergodic.

\begin{theorem}[\cite{BezuglyiKarpelKwiatkowski2018}]\label{uniq_erg} 
A Bratteli diagram $B = (V,E)$ is uniquely ergodic if and only if there exists a telescoping $B'$ of $B$ such that
\begin{equation}\label{2.3b}
\lim_{n\to \infty}\ \max_{v, v' \in V_{n+1}} \left(\sum_{w \in V_n}
\left|f_{vw}^{(n)} - f_{v'w}^{(n)}\right| \right) = 0, 
\end{equation}
where $f_{vw}^{(n)}$ are the entries of
 the stochastic matrix $F_n$ defined by the diagram $B'$.
\end{theorem}

The proof is based on the representation of an invariant measure as a point of 
the inverse limit of the
sequence $(F_n^T, \Delta^{(n)}_{\infty})$ (see Section~\ref{Section_Inv_measures}). We use the fact that $B$ is uniquely ergodic if and only if the set $\Delta_{\infty}^{(n)}$ is a singleton for all $n = 1,2,\ldots$  and that the polytope $\Delta_m^{(n)}$ is the convex hull of the vectors $\{\ov g^{(n+m,n)}(v)\}_{v \in V_{n+m+1}}$ for all $m \in \mathbb{N}$.

The following statement is a corollary of Theorem~\ref{uniq_erg} and provides a sufficient condition for a Bratteli diagram to be uniquely ergodic. Note that this condition does not require telescoping.

\begin{theorem}[\cite{BezuglyiKarpelKwiatkowski2018}]\label{suff_cond_ue}
Let $B$ be a Bratteli diagram of arbitrary rank with stochastic incidence matrices $F_n$ and let
 $$
 m_n = \min_{v \in V_{n+1}, w \in V_n} f_{vw}^{(n)}.
 $$
 If  
 \begin{equation}\label{eq_suff_cond_ue}
 \sum_{n = 1}^{\infty} m_n = \infty,
\end{equation}
 then $B$ is uniquely ergodic.
\end{theorem}

A number of sufficient conditions for unique ergodicity of a finite rank Bratteli-Vershik system are obtained in~\cite{BezuglyiKwiatkowskiMedynetsSolomyak2013}. Here we present some of them.

\begin{definition}(see e.g. \cite{Hartfiel2002})
(i) For two positive vectors $x,y\in\mathbb R^d$, the {\it projective metric} is defined by the formula
$$
D(x,y) = \ln\max_{i,j}\frac{x_iy_j}{x_jy_i} = \ln \frac{\max_i\frac{x_i}{y_i}}{\min_j \frac{x_j}{y_j}},
$$
where $(x_i)$ and $(y_i)$ are entries of the vectors $x$ and $y$.

(ii) For a non-negative matrix $A$, the {\it Birkhoff contraction
coefficient} is $$\tau(A) = \sup_{x,y>0}\frac{D(Ax,Ay)}{D(x,y)}.$$
\end{definition}

\begin{theorem}[\cite{BezuglyiKwiatkowskiMedynetsSolomyak2013}]\label{theoremUniqueErgodicityInTermsTau}
Let $B$ be a simple Bratteli diagram of finite rank
with incidence matrices $\{\tl F_n\}_{n\geq 1}$. Let $\tl A_n = \tl F_n^T$.
Then the diagram $B$ is uniquely ergodic if and only if
$$\lim_{n\to
\infty}\tau(\tl A_m\ldots \tl A_n) = 0\mbox{ for every }m.
$$
\end{theorem}

For a positive matrix $A = (a_{i,j})$, let
$$\phi(A) = \min_{i,j,r,s} \frac{a_{i,j}a_{r,s}}{a_{r,j}a_{i,s}}.
$$
If $A$ has a zero entry, then, by definition, we set $\phi(A)=0$.
As noticed in \cite{Hartfiel2002},
 $$\tau(A) = \frac{1-\sqrt{\phi(A)}}{1+\sqrt{\phi(A)}}$$
when  $A$ has a nonzero entry in each row.

The following result gives sufficient conditions for unique ergodicity which can be easily verified when a diagram is given by a sequence of incidence matrices.

\begin{proposition}[\cite{BezuglyiKwiatkowskiMedynetsSolomyak2013}]\label{PropositionSufficientConditionsUniqueErgodicity}
Let $\{\tl A_n\}_{n\geq 1} = \tl F^T_n$ be primitive incidence matrices of a finite rank diagram $B$.

(1) If $$\sum_{n = 1}^\infty \sqrt{\phi(A_n)} = \infty,$$
then $B$ admits a unique invariant probability measure.

(2) If
$$
\sum_{n = 1}^\infty \left(\frac{\tl m_n}{\tl M_n} \right)=\infty,
$$
where $\tl m_n$ and $\tl M_n$ are the smallest and the largest entry of
$\tl A_n$ respectively, then $B$ admits a unique invariant probability
measure.

(3) Let $||A||_1 := \sum_{i,j}|a_{i,j}|$.  If  $||\tl F_n||_1\leq Cn$ for some $C>0$ and all sufficiently
large $n$, then the diagram admits a unique invariant probability
measure. In particular, this result holds if the diagram has only finitely
many different incidence matrices.

\end{proposition}

\subsection{Examples}

The following examples illustrate the results of Subsection~\ref{subsect_finrank_ue}. In particular, Examples~\ref{ex 3.2} and \ref{non-simple-stat-ex} show that telescoping and using the stochastic incidence matrix are crucial for Theorem~\ref{uniq_erg}. 

\begin{example} \label{ex 3.2} Let $B$ be a Bratteli diagram with incidence
  matrices
$$
\tl F_n =
\begin{pmatrix}
n & 1\\
1 & n
\end{pmatrix}, \qquad n \in \mathbb N.
$$
Then the diagram $B$ is uniquely ergodic, see details in
 \cite{BezuglyiKwiatkowskiMedynetsSolomyak2013},
 \cite[Example 3.6]{AdamskaBezuglyiKarpelKwiatkowski2016},  and
 \cite{FerencziFisherTalet2009}.

 
Notice that $B$ has the ERS property (see Example~\ref{ExampleERS1}). Hence the
 corresponding stochastic incidence matrices are:
$$
F_n =
\begin{pmatrix}
1 - \dfrac{1}{n+1} & \dfrac{1}{n+1}\\
\\
\dfrac{1}{n+1} & 1 - \dfrac{1}{n+1}
\end{pmatrix}.
$$
Obviously, without telescoping, for $B$ the limit in~(\ref{2.3b})
 equals $2$. However, the telescoping procedure reveals that the diagram is in fact uniquely ergodic.

Suppose we have an ERS diagram with $2 \times 2$ stochastic incidence
 matrices
$$
F'_n =
\begin{pmatrix}
a_n & b_n\\
b_n & a_n
\end{pmatrix}.
$$
As before, let $G_{(n,n+m)} = (g_{vw}^{(n+m,n)})$ be the corresponding
 product matrix.
It can be easily proved  by induction that,  for arbitrary
$n, m  \in \mathbb{N}$, the following formula holds:
$$
S(n,m) = \sum_{w \in V_n} \left|g_{vw}^{(n+m,n)} - g_{v'w}^{(n+m,n)}\right| =
 2\prod_{i = 0}^{m}\left|(a_{n+i} - b_{n+i})\right|.
$$

In the case of the diagram $B$, we obtain
$$
S(n,m) = 2\prod_{i = 0}^{m}\left(1 - \frac{2}{n+i+1}\right), \quad
n, m \in \mathbb N.
$$
Since the series $\sum_{n=1}^{\infty} n^{-1}$ diverges, we see
that
$$
S(n,m) \to 0, \qquad m \to \infty.
$$
 Choose a decreasing sequence $(\varepsilon_k)$ such that
 $\varepsilon_k \to 0$ as $k \to \infty$. For $n = n_1$ and
 $\varepsilon_1$, find $m_1$ such that $S(n_1, m_1) < \varepsilon_1$.
 Set $n_2 = n_1 + m_1$. For  $\varepsilon_2$, find $m_2$ such that
 $S(n_2, m_2) < \varepsilon_2$. Set $n_3 = n_2 + m_2$. Continuing
 in the same manner, we construct a sequence $(n_k)$ such
 that $S(n_k, n_{k+1} - n_k) < \varepsilon_k$. Telescope the diagram
 with respect to the levels $(n_k)$. By Theorem
 \ref{uniq_erg}, we conclude that the diagram $B$ is uniquely ergodic. Notice that the diagram $B$ also satisfies the sufficient condition of unique ergodicity~(\ref{eq_suff_cond_ue}).


\end{example}

\begin{example}\label{Ex_generalized}
 Let $B$ be a simple Bratteli diagram with incidence matrices
$$
\tl F_n = \left(
  \begin{array}{cccc}
    f_1^{(n)} & 1 & \cdots & 1 \\
    1 & f_2^{(n)} & \cdots & 1 \\
    \vdots & \vdots & \ddots & \vdots \\
    1 & 1 & \cdots & f_d^{(n)} \\
      \end{array}
\right).
$$
Let $q_n = \mbox{max}\{\tl f_i^{(n)}\tl f_j^{(n)} : i\neq j \}$. By Proposition~\ref{PropositionSufficientConditionsUniqueErgodicity}, if for $\tl A_n = \tl F_n^T$
$$
\sum_{n=1}^\infty \sqrt{\phi(\tl A_n)} = \sum_{n=1}^\infty
\frac{1}{\sqrt q_n} =\infty,
$$
then there is a unique invariant probability measure on $B$.
\end{example}

\begin{example}\label{non-simple-stat-ex}
Let $B$ be the stationary non-simple Bratteli diagram defined by the incidence
  matrices
$$
\tl F_n = \tl F =
\begin{pmatrix}
3 & 0\\
1 & 2
\end{pmatrix}
$$
for every $n \in \mathbb N$. It is well known that $B$ has a unique finite
 ergodic  measure supported by the 3-odometer (see e.g.
 \cite{BezuglyiKwiatkowskiMedynetsSolomyak2010}). We show that $B$ 
 satisfies the condition of unique ergodicity formulated in 
 Theorem \ref{uniq_erg}. It is easy to check that the $n$-th power of 
 $\tl F$ is
$$
\tl F^{n} =  
\begin{pmatrix}
3^n & 0\\
3^n - 2^n & 2^n
\end{pmatrix}.
$$
Hence the entries of the matrix $\tl F$ do not satisfy (\ref{2.3b}) even after 
taking products of these matrices (which corresponds to telescoping of $B$).
 Notice that $B$ has the 
ERS property. For any $n \in \mathbb{N}$ and a vertex $w \in V_n$, we 
have $h_w^{(n)} = 3^n$. Therefore, the corresponding stochastic incidence
 matrix and its $n$-th power are
$$
F =
\begin{pmatrix}
1 & 0\\
\frac{1}{3} & \frac{2}{3}
\end{pmatrix}
$$
and
$$
F^{n} =  
\begin{pmatrix}
1 & 0\\
1 - \dfrac{2^n}{3^n} & \dfrac{2^n}{3^n}
\end{pmatrix}.
$$
Hence, we see that $B$ satisfies (\ref{2.3b}) and is uniquely ergodic.
Note that $B$ does not satisfy the sufficient condition of unique ergodicity~(\ref{eq_suff_cond_ue}) and that Theorem~\ref{theoremUniqueErgodicityInTermsTau} and Proposition~\ref{PropositionSufficientConditionsUniqueErgodicity} are not applicable since $B$ is not simple.
\end{example}

\section{Finitely ergodic Cantor dynamical systems}\label{finrank}

This section is mostly devoted to the study of aperiodic Cantor dynamical systems which can be represented by Bratteli diagrams with uniformly bounded number of vertices on each level. It is an open question which classes of Cantor dynamical systems admit such a representation.

\subsection{Finitely ergodic subshifts}\label{subs FE subshifts}

As in Section \ref{Section_Uniq_Erg}, we begin with the case of finitely
ergodic (minimal) subshifts.  The recent progress made in \cite{CyrKra2019},
\cite{DamronFickenscher2017}, \cite{DamronFickenscher2019} 
essentially  improved the known results on the bounds of the cardinality
of the set $E(X, S)$ of ergodic invariant measures.

In \cite{Boshernitzan1984}, the following remarkable results were proved. 
\begin{theorem}[\cite{Boshernitzan1984}]  Let $(X, S)$ be a minimal subshift on a finite alphabet
$\A$. 
\\
(i) If 
$$
\liminf_{n\to \infty} \frac{p_X(n)}{n} = \alpha,
$$
then $|E(X, S)| \leq [\alpha]$, where $[\alpha]$ the integer part of $\alpha$.
\\
(ii) If 
$$
\limsup_{n\to \infty} \frac{p_X(n)}{n} = \alpha
$$ 
and $\alpha \geq 2$, then $|E(X, S)| \leq [\alpha] - 1$.
\end{theorem}

The  results given in  \cite{FerencziMonteil2010} extend Boshernitzan's 
bounds to the so called $K$-deconnectable symbolic systems. We cite only 
one of these results from that paper here.

\begin{theorem}[\cite{Monteil2009}]
Let $K \geq 3$ be an integer. A minimal symbolic system $(X, S)$ such that
$$
\limsup_{n\to \infty} \frac{p_X(n)}{n} < K
$$
admits at most $K-2$ ergodic invariant measures.
\end{theorem}

In \cite{DamronFickenscher2017}, the authors continued this line of study
of the set $E(X, S)$ and considered complexity functions with eventually
constant growth condition. By definition, this means that the 
complexity function $p_X(n)$ of a minimal subshift satisfies the
condition: for some $K\in \N$ and all $n \geq n_0$
\begin{equation}\label{eq ev const}
p_X(n+1) - p_X(n) = K.
\end{equation}
Equivalently, $p_X(n) = Kn + C$ for all $n \geq n_0$ where a constant 
$C \in \N_0$.

\begin{theorem} [\cite{DamronFickenscher2017}]
 If the complexity function of a minimal subshift $(X, S)$
satisfies eventually constant growth condition \ref{eq ev const} with
$K\geq 4$, then $|E(X, S)| \leq K - 2$.
\end{theorem}

In the very recent paper \cite{DamronFickenscher2019}, the authors 
addressed the old \textit{question} asked by Boshernitzan. 
Let $(X, T)$ be a minimal interval exchange transformation (IET) defined by a
permutation of $d$
subintervals. Due to Katok \cite{Katok1973} and Veech \cite{Veech1978}, 
it is known that 
$$
|E(X, T)| \leq \left[\frac{d}{2}\right].
$$  
Can the bound $\frac{d}{2}$ for the IET $|E(X, T)|$  be shown 
combinatorially  using  a symbolic realization $(Y, S)$ of $(X, T)$?

Following \cite{DamronFickenscher2019}, let us make the following 
assumption on the language $\La_X$ of a minimal subshift $(X, S)$. 
A word $w \in  \La_X$ is left special if there are distinct letters
$a, a' \in  \A$ such that $aw$ and $a'w$ belong to $\La_X$. Likewise, $w$
 is right special if $wb$ and $wb'$ exist in the language $\La_X$ for distinct
  letters $b, b'$. A word $w$ is bispecial if it
is both left and right special. A bispecial word is  called regular bispecial if 
only one left extension of $w$ is right special and only one right extension 
of $w$ is left special. The language $\La_X$
(or equivalently$ (X, S)$) satisfies the regular bispecial condition  if all large
enough bispecial words are regular. The regular bispecial condition 
 implies the constant growth condition above for some $K$.
All subshifts that arise from interval exchanges satisfy this property, see
\cite{FerencziZamboni2008}.

The following main result from \cite{DamronFickenscher2019} is motivated 
by the Boshernitzan's question.
\begin{theorem} [\cite{DamronFickenscher2019}]
Let $(X, S)$ be a transitive subshift satisfying the regular bispecial condition
with growth constant $K$. Then 
$$
|E(X, S)| \leq  \frac{K+1}{2}.
$$ 
\end{theorem}

We finish this subsection by pointing out an  interesting application of  
complexity functions. It turns out that by means of the complexity function 
one can also estimate
the number of generic measures (which are not necessarily ergodic ones). 
We follow here the paper \cite{CyrKra2019}. We remark that the considered
subshifts are not assumed to be minimal.

\begin{theorem} [\cite{CyrKra2019}] (1) Let $(X, S)$ be a subshift such 
that 
$$
\liminf_{n\to \infty} \frac{p_X(n)}{n } < K
$$
for some integer $K$. Then $(X, S)$ has at most $K-1$ distinct, non-atomic,
generic measures.
\\
(2) Suppose $(X, S)$ is a subshift satisfying the condition
$$
\limsup_{n\to \infty} \frac{p_X(n)}{n } < K
$$
for some integer $K$.  If $(X, S)$ has a generic measure $\mu$ and a generic
point $x_\mu$ for which the orbit closure $(\overline{Orb_S(x_\mu)}, S)$
is not uniquely ergodic, then $(X, S)$ has at most $K-2$ distinct, non-atomic,
generic measures.
\end{theorem}

\subsection{Stationary Bratteli diagrams}\label{stationary}

In this subsection, we give an explicit description of all ergodic probability invariant measures on stationary Bratteli diagrams. Note that the class of minimal homeomorphisms which can be represented by stationary Bratteli diagrams is constituted by minimal substitution dynamical systems and odometers~\cite{Forrest1997, DurandHostSkau1999}. In \cite{BezuglyiKwiatkowskiMedynets2009}, the analogue of the above mentioned result was proved for aperiodic homeomorphisms.

The paper~\cite{BezuglyiKwiatkowskiMedynetsSolomyak2010} contains an explicit
description of all ergodic invariant probability measures on a stationary Bratteli diagram $B$.
Let $\tl F = (\tl f_{vw})_{v,w \in V}$ be the $K \times K$ incidence matrix of the diagram $B$.
Identify the set of vertices $V_n$ on each level $n \geq 1$ with $\{1, \ldots, K\}$.
In this subsection, by $\ov x$ we denote a  vector, either column or row one, it will
be either mentioned explicitly, or understood from the context.

The incidence matrix $\tl F$ defines a directed
graph $G(\tl F)$: the set of the vertices of $G(\tl F)$ is
equal to $\{1, \ldots, K\}$ and there is a directed edge from a vertex $v$ to a vertex $w$
if and only if $\tl f_{vw} > 0$. The vertices $v$ and $w$ are equivalent (we
write $v \sim w$) if either $v = w$ or there is a path in $G(\tl F)$ from $v$
to $w$ and also a path from $w$ to $v$. Let $\mathcal{E}_1,\ldots,
\mathcal{E}_m$ denote all equivalence classes in $G(\tl F)$. We will also
identify $\mathcal{E}_{\alpha}$ with the corresponding subsets of $V$. We
write $\mathcal{E}_{\alpha} \succeq \mathcal{E}_{\beta}$ if either
$\mathcal{E}_{\alpha} = \mathcal{E}_{\beta}$ or there is a path in
$G(\tl F)$ from a vertex of $\mathcal{E}_{\alpha}$ to a vertex of
$\mathcal{E}_{\beta}$. We write $\mathcal{E}_{\alpha} \succ
\mathcal{E}_{\beta}$ if $\mathcal{E}_{\alpha} \succeq
\mathcal{E}_{\beta}$ and $\mathcal{E}_{\alpha} \neq
\mathcal{E}_{\beta}$. Every class $\mathcal{E}_{\alpha}$, $\alpha = 1,
\ldots,m$, defines an irreducible submatrix $\tl F_{\alpha}$ of $\tl F$
obtained by restricting $\tl F$ to the set of vertices from
$\mathcal{E}_{\alpha}$. Let $\rho_{\alpha}$ be the spectral radius of
$\tl F_{\alpha}$, i.e.
$$
\rho_{\alpha} = \max\{|\lambda|: \lambda \in
\mathsf{spec}(\tl F_{\alpha})\},
$$
where by $\mathsf{spec}(\tl F_{\alpha})$ we mean the set of all complex
numbers $\lambda$ such that there exists a non-zero vector $\ov x =
(x_v)_{v \in \mathcal{E}_{\alpha}}$ satisfying $\ov x \tl{F}_{\alpha} =
\lambda \overline x$.

A class $\mathcal{E}_{\alpha}$ is called \textit{distinguished} if
\begin{equation}\label{eq disting class}
\rho_{\alpha} > \rho_{\beta} \ \ \mathrm{whenever} \ \ \
\mathcal{E}_{\alpha} \succ \mathcal{E}_{\beta}
\end{equation}
(in~\cite{BezuglyiKwiatkowskiMedynetsSolomyak2010} the notion of being
distinguished is defined in an opposite way because it is based on the matrix
transpose to the incidence matrix).

The real number $\lambda$ is called a \textit{distinguished eigenvalue} if
 there exists a non-negative left-eigenvector $\ov x = (x_v) \in \mathbb{R}^K$ such
 that $\ov x \tl F = \lambda \ov x$. It is known (Frobenius theorem) that $
 \lambda$ is a distinguished eigenvalue if and only if $\lambda =
 \rho_{\alpha}$ for some distinguished class $\mathcal{E}_{\alpha}$.
 Moreover, there is a unique (up to scaling) non-negative eigenvector
 $\ov  x(\alpha) = (x_v)_{v \in V}$, $\ov x(\alpha)\tl F = \rho_{\alpha}
 \ov  x(\alpha)$ such that $x_v > 0$ if and only if there is a path from a
  vertex  of $\mathcal{E}_\alpha$ to the vertex $v$. The distinguished class
   $\alpha$ defines a measure $\mu_{\alpha}$ on $B = (V,E)$ as follows:
$$
\mu_{\alpha}(X_v^{(n)}) = \frac{x_v}{\rho_{\alpha}^{n-1}} h_v^{(n)}, \ \
\ v \in V_n = V.
$$

\begin{theorem}[\cite{BezuglyiKwiatkowskiMedynetsSolomyak2010}]
Let $B$ and $\{\mu_{\alpha}\}$ be as above, where $\alpha$ runs over all distinguished vertex classes. Then the measures $\{\mu_{\alpha}\}$ are exactly all probability ergodic $\mathcal{R}$-invariant
measures for the stationary Bratteli diagram $B$.
\end{theorem}

For instance, in Example~\ref{non-simple-stat-ex}, there is only one distinguished class of vertices which corresponds to the first vertex of the diagram on each level.

\begin{remark} In \cite{BezuglyiKwiatkowskiMedynetsSolomyak2010} it was shown that non-distinguished vertex classes correspond exactly to infinite ergodic invariant measures which are finite on at least one open set.
\end{remark}
 

\subsection{Finite rank Bratteli diagrams}

In this subsection, we give the necessary and sufficient conditions to determine the exact number of probability ergodic invariant measures on Bratteli diagrams of finite rank and describe the supports of these measures. 

\begin{definition}
A Cantor dynamical system $(X,S)$ has the {\it topological rank} $K>0$ if it admits a Bratteli-Vershik model $(X_B,\varphi_B)$ such that the number of vertices of the diagram $B$ at each level $V_n, n \geq 1$ is not greater than $K$ and $K$ is the least possible number of vertices for any Bratteli-Vershik realization.
\end{definition}

If a system $(X,S)$ has the rank $K$, then, by an appropriate telescoping, we can assume that the diagram $B$ has exactly $K$ vertices at each level. 

In~\cite{BezuglyiKwiatkowskiMedynetsSolomyak2013}, the structure of invariant measures on finite rank Bratteli diagrams is considered. In particular, it is shown that every ergodic invariant measure (finite or ``regular'' infinite) can be obtained as an extension from a simple uniquely ergodic vertex subdiagram.
Everywhere below the term ``measure'' stands for an $\mathcal R$-invariant measure. By an infinite measure we mean any $\sigma$-finite non-atomic measure which is finite (non-zero) on some clopen set. The support of each ergodic measure for a Bratteli diagram of finite rank turns out to be the set of all paths that stabilize in some subdiagram, which geometrically can be seen as a ``vertical'' subdiagram, i.e. the paths will eventually stay in the subdiagram. Furthermore, these subdiagrams are pairwise disjoint for different ergodic measures. It is shown in~\cite{BezuglyiKwiatkowskiMedynetsSolomyak2013}, that for any finite rank diagram $B$ one can find finitely many vertex subdiagrams $B_\alpha$ such that each finite ergodic measure on $X_{B_\alpha}$ extends to a (finite or infinite) ergodic measure on $X_B$. It is also proved that each ergodic measure (both finite and infinite) on $X_B$ is obtained as an extension of a finite ergodic measure from some $X_{B_\alpha}$. The following theorem holds:

\begin{theorem}[\cite{BezuglyiKwiatkowskiMedynetsSolomyak2013}]\label{TheoremGeneralStructureOfMeasures}
Let $B$ be a Bratteli diagram of finite rank $K$.
The diagram $B$ can be telescoped in such a way that for every
probability ergodic measure $\mu$ there exists a subset $W_\mu$ of
vertices from $\{1,\ldots,K\}$ such that the support of $\mu$ consists of
all infinite paths that eventually go along the vertices of $W_\mu$
only. Furthermore,

(i) $W_\mu\cap W_\nu = \emptyset$ for different ergodic measures
$\mu$ and $\nu$;

(ii) given a probability ergodic measure $\mu$, there exists a constant $\delta>0$
such that for any $w\in W_\mu$ and any level $n$
$$\mu(X_w^{(n)})\geq \delta;
$$

(iii) the subdiagram generated by $W_\mu$ is simple and uniquely
ergodic. The only ergodic measure on the path space of the subdiagram
is the restriction of measure $\mu$;

(iv) if a probability ergodic measure $\mu$ is the extension of a measure
 from the vertical subdiagram determined by a proper subset $W\subset\{1,\ldots,K\}$, then
$$\lim_{n\to\infty}\mu(X_w^{(n)}) = 0\mbox{ for all }w\notin W.$$
\end{theorem}

Condition (ii) can be used in practice to determine the support
of an ergodic measure $\mu$.

\medskip
The following theorem plays an important role in the study of ergodic
measures and their supports. For a finite rank Bratteli diagram $B$, it
 describes how extreme points of $\Delta_{\infty}^{(1)}$ determine 
 subdiagrams of $B$.

\begin{theorem}[\cite{BezuglyiKarpelKwiatkowski2018}]\label{prop subdiagrams}
Let $B$ be a Bratteli diagram of rank $K$, and let $B$ have $l$ probability
ergodic invariant measures, $1 \leq l \leq K$. Let $\{\ov y_1, ... , \ov y_l\}$ denote the extreme vectors in $\Delta_{\infty}^{(1)}$.
Then, after telescoping and
renumbering  vertices, there exist exactly $l$ disjoint subdiagrams $B_i$
(they share no vertices other than the root) with the corresponding sets of vertices $\{V_{n,i}\}_{n = 0}^{\infty}$
such that

(a) for every $i = 1,...,l$ and any $n,m >0$, $|V_{n,i}| = |V_{m,i}| > 0$,
while the set $V_{n,0} = V_n \setminus \bigsqcup\limits_{i = 1}^l V_{n,i}$ may be, in particular, empty;

(b) for any $i = 1, \ldots, l$ and any choice of $v_n \in V_{n,i}$, the extreme vectors $\ov y^{(n)}
(v_n) \in \Delta_n^{(1)}$ converge to the extreme vector $\ov y_i \in
\Delta_{\infty}^{(1)}$.

In general, the diagram $B$ can have up to $K - l$ disjoint subdiagrams
$B'_j$ with vertices $\{V'_{n,j}\}_{n = 0}^{\infty}$ such that they are
 also disjoint with subdiagrams $B_i$ and for any $w_n \in V'_{n,j}$, the
  extreme vectors $\ov y^{(n)}(w_n) \in \Delta_n^{(1)}$ converge to a
  non-extreme vector $\ov z \in \Delta_{\infty}^{(1)}$.
\end{theorem}

For a finite rank Bratteli diagram, one can describe subdiagrams that support ergodic measures in terms of the stochastic incidence matrices of the diagram. For the next theorem,  we will need the following definition and notation.

\begin{definition}\label{def vanishing blocks}
For a Bratteli diagram $B$,  we say that a sequence of proper subsets
$U_n \subset V_n$ defines \textit{blocks of vanishing weights
(or vanishing blocks)} in the stochastic incidence matrices $F_n$ if
$$
 \sum_{w \in U^c_n, v \in U_{n+1}} f_{vw}^{(n)} \rightarrow 0, \qquad
 n \rightarrow \infty
 $$
 where $U_n^c = V_n \setminus U_n$.

If additionally, for every sequence of vanishing blocks $(U_n)$, there exists
a constant  $0 < C_1 < 1$ such that, for sufficiently large $n$,
\begin{equation}\label{eq vanish blocks}
\min_{v\in U_{n+1}} \sum_{u \in U_{n}^c} f_{vu}^{(n)}
\geq C_1 \max_{v\in U_{n+1}} \sum_{u \in U_{n}^c}  f_{vu}^{(n)},
\end{equation}
then we say that the stochastic incidence matrices $F_n$ of $B$ have
 \textit{regularly vanishing blocks}.
\end{definition}

Set
$$
\ov a_j^{(n)} = \frac{1}{|V_{n,j}|} \sum_{w \in V_{n,j}} \ov y^{(n)}(w),
\quad j = 0, 1, \ldots, l,
$$
where the subsets $V_{n,j}$ are defined as in Theorem
\ref{prop subdiagrams}.
Then $\ov a_j^{(n)} \in \Delta_{n,j}^{(1)} :=
\textrm{Conv}\{\ov y^{(n)}(w),
w \in V_{n,j}\}$, the convex hull of the set $\{\ov y^{(n)}(w),
w \in V_{n,j}\}$.   The sets $\Delta_{n,j}^{(1)}$ are subsimplices of
$\Delta_n^{(1)}$, $j = 0, 1, \ldots, l$.
We observe that
 $$
 \max_{\ov a \in \Delta_{n,0}^{(1)}}
\textrm{dist}(\ov a, \Delta_{\infty}^{(1)}) \rightarrow 0
$$
as $n \rightarrow \infty$, where $\Delta_{\infty}^{(1)} =
\bigcap_{n=1}^{\infty}\Delta_n^{(1)}$.
\medskip


In the theorem below, we assume that stochastic incidence matrices of a Bratteli diagram have the property of {\it regularly vanishing blocks} and apply it to the case when $U_n = V_{n,i}$ for some $i = 1, \ldots, l$. The blocks of the matrices corresponding to the edges that connect vertices from outside of the supporting subdiagram $B_i$ to the vertices of $B_i$ are the blocks of vanishing weights. Note, that the second part of the theorem does not require the stochastic incidence matrices of a Bratteli diagram to have the property of regularly vanishing blocks.

\begin{theorem}[\cite{BezuglyiKarpelKwiatkowski2018}]\label{main1} Let $B$ be a Bratteli diagram of rank $K$ such
that the incidence matrices $F_n$ have the property of regularly vanishing
blocks. If $B$ has exactly
$l\ (1 \leq l \leq K)$ ergodic invariant probability measures, then, after
 telescoping, the set $V_n$ can be partitioned into subsets $\{V_{n,1},
 \ldots,V_{n,l},V_{n,0}\}$  such that

\smallskip
(a) $V_{n,i} \neq \emptyset$ for $i = 1,\ldots,l$;

\smallskip
(b) $|V_{n,i}|$ does not depend on $n$, i.e.,
$|V_{n,i}| = k_i$ for $i = 0,1, \ldots,l$ and $n \geq 1$;

\smallskip
(c) for $j = 1,\ldots,l$,
$$
\sum_{n = 1}^\infty \left(1 - \min_{v \in V_{n+1,j}}\sum_{w \in V_{n,j}}
f_{vw}^{(n)}\right)< \infty;
$$

\smallskip
(d)  for $j = 1,\ldots,l$,
$$
\max_{v,v' \in V_{n+1,j}} \sum_{w \in V_{n}} \left|f_{vw}^{(n)} -
f_{v'w}^{(n)}\right| \rightarrow 0
$$
as $n \rightarrow \infty$;

\smallskip

\ \ (e1)
for every $w \in V_{n,0}$
$$
\mathrm{vol}_{l} S(\ov a_1^{(n)}, \ldots, \ov a_l^{(n)}, \ov y^{(n)}(w))
 \rightarrow 0
$$
as $n \rightarrow \infty$, where $S$ is a simplex with extreme points
$\ov a_1^{(n)}, \ldots, \ov a_l^{(n)}, \ov y^{(n)}(w)$, and
$\mathrm{vol}_l(S)$ stands for the volume of $S$;

 (e2) for every $v \in V_{n+1,0}$ and  for sufficiently large $n$,
there exists some $C>0$ such that, for every $j = 1, \ldots,l$,
$$F_v^{(n,j)} =  \sum_{w  \in V_{n,j}}  f_{vw}^{(n)} < 1 -C.
$$

\medskip
Conversely, let $B$ be a Bratteli diagram of finite rank $K \geq 2$ with nonsingular 
stochastic incidence matrices $(F_n)$. Suppose that  after telescoping $B$ 
satisfies conditions $(a) - (e2)$. Then $B$ has $l$ ergodic probability invariant measures.
\end{theorem}

\begin{remark}
Condition (d) of Theorem~\ref{main1} guarantees that the subdiagrams $B_i$, $i = 1, \ldots, l$ corresponding to the vertices from $V_{n,i}$, are uniquely ergodic, but it does not guarantee that the subdiagrams $B_i$ are simple. One can reduce
  $B_i$ to the smallest possible simple and uniquely 
ergodic subdiagrams such that the obtained subdiagrams are the same 
as considered in Theorem~\ref{TheoremGeneralStructureOfMeasures}. For instance, in 
Example~\ref{non-simple-stat-ex}, one can take $V_{n,1} = V_n$ and 
$V_{n,0} = \emptyset$ for all $n$. After reduction, we obtain that the new set 
$V'_{n,1}$ consists only of the first vertex on each level $n$, and $V'_{n,0}$ 
consists of the second one. Condition (c) of Theorem~\ref{main1} yields that for every $i = 1, \ldots, l$, the extension of the unique invariant measure $\mu_i$ on $B_i$ to the measure $\wh \mu_i$ on $B$ is finite. Conditions $(e1)$ and $(e2)$ guarantee that there are no more finite ergodic invariant measures on $B$ except for $\wh \mu_1, \ldots \wh \mu_l$.
\end{remark}

The following theorem gives a criterion for the existence of $K$ probability ergodic invariant measures on a Bratteli diagram of rank $K$. This criterion was proved in \cite{AdamskaBezuglyiKarpelKwiatkowski2016} for the case of Bratteli diagrams with ERS property, but actually it can be reproved in terms of stochastic incidence matrices $(F_n)$ without the ERS property requirement. 
\begin{theorem}\label{Thm_max_numb_meas}
Let $B = (V,E)$ be a Bratteli diagram of rank $K \geq 2$; identify $V_n$ with $\{1, ... ,K\}$ for any $n \geq 1$. Let $F_n = (f_{i,j}^{(n)})$ form a sequence of stochastic incidence matrices of $B$. Suppose that  $\rank\ F_n = K$ for all $n$. 
Denote
$$
z^{(n)} =
\det
\begin{pmatrix}
f_{1,1}^{(n)} & \ldots & f_{1,{k}}^{(n)}\\
\vdots & \ddots & \vdots \\
f_{{k},1}^{(n)} & \ldots & f_{{k},{k}}^{(n)}
\end{pmatrix}.
$$
Then there exist exactly $K$ ergodic invariant measures on $B$ if and only if
$$
\prod_{n = 1}^{\infty} |z^{(n)}| > 0,
$$
or, equivalently,
$$
\sum_{n = 1}^{\infty}(1 - |z^{(n)}|) < \infty.
$$
\end{theorem}

\subsection{Examples}

\begin{example}[Stationary Bratteli diagrams]
This example illustrates Theorem~\ref{main1}.
For stationary Bratteli diagrams (see Subsection~\ref{stationary}), we relate the distinguished classes of vertices to the subsets $V_{n,j}$ mentioned in Theorem~\ref{main1}. 


\begin{proposition}[\cite{BezuglyiKarpelKwiatkowski2018}] Let $B = (V,E)$ be a stationary Bratteli diagram and $V_{n,j}$, $j = 1,\ldots,l$ be subsets of vertices defined in Theorem~\ref{main1}. Then the distinguished classes $\alpha$ (as subsets of $V$)
coincide with the sets $V_{n,j}$, $j = 1,\ldots,l$.
\end{proposition}

The proof of the above proposition uses the representation of the incidence matrix $\tl F$ in the Frobenius normal form
(similarly to the way it was done in~
\cite{BezuglyiKwiatkowskiMedynetsSolomyak2010}):
$$
F =\left(
  \begin{array}{ccccccc}
    \tl F_1 & 0 & \cdots & 0 & 0 & \cdots & 0 \\
    0 & \tl F_2 & \cdots & 0 & 0 & \cdots & 0 \\
    \vdots & \vdots & \ddots & \vdots & \vdots & \cdots& \vdots \\
    0 & 0 & \cdots & \tl F_s & 0 & \cdots & 0 \\
    Y_{s+1,1} & Y_{s+1,2} & \cdots & Y_{s+1,s} & \tl F_{s+1} & \cdots & 0
    \\
    \vdots & \vdots & \cdots & \vdots & \vdots & \ddots & \vdots \\
    Y_{m,1} & Y_{m,2} & \cdots & Y_{m,s} & Y_{m,s+1} & \cdots & \tl F_m
    \\
  \end{array}
\right),
$$
where all $\{\tl F_i\}_{i = 1}^m$ are irreducible square matrices, and for any
$j = s+1,\ldots,m$, at least one of the matrices $Y_{j,u}$ is non-zero. All
classes $\{\mathcal{E}_{\alpha}\}_{\alpha = 1}^s$ $(s \geq 1)$, are
distinguished (there is no $\beta$ such that $\alpha > \beta$). For every
$\alpha \geq s+1$ such that $\mathcal{E}_{\alpha}$ is a distinguished class
and for every $1 \leq \beta < \alpha $ we have either $\mathcal{E}_{\beta}
\prec \mathcal{E}_{\alpha}$ and $\rho_{\beta} < \rho_{\alpha}$, or there
is no relation between $\mathcal{E}_{\alpha}$ and $\mathcal{E}_{\beta}$.
Then the Perron-Frobenius theorem is used to show that the sets $V_{n,j}$, $j = 1,\ldots,l$, coincide with the distinguished classes $\mathcal{E}_{\alpha}$.

\end{example}

\begin{example}
Let $B$ be a Bratteli diagram of rank $2$ with incidence matrices
$$
\tl F_n =
\begin{pmatrix}
n^2 & 1\\
1 & n^2
\end{pmatrix}, \qquad n \in \mathbb N.
$$
Since $B_2$ has the ERS property (see Example~\ref{ExampleERS1}), the corresponding stochastic incidence matrices are:
$$
F_n =
\begin{pmatrix}
1 - \dfrac{1}{n^2+1} & \dfrac{1}{n^2+1}\\
\\
\dfrac{1}{n^2+1} & 1 - \dfrac{1}{n^2+1}
\end{pmatrix}.
$$
Since the series
$$
\sum_{n=1}^{\infty} \dfrac{1}{n^2+1}
$$
converges, by Theorem~\ref{main1} the diagram $B$  has two probability ergodic invariant measures.
It also is easy to see that the diagram satisfies the condition of Theorem~\ref{Thm_max_numb_meas}  (see also Proposition 3.1 in \cite{AdamskaBezuglyiKarpelKwiatkowski2016}). 

This example can be generalized to the case of Bratteli diagrams of rank $K \geq 2$ by using Example~\ref{Ex_generalized} and choosing the appropriate values for $f_i^{(n)}$, $i = 1,\ldots,n$. 
\end{example}

\section{Infinite rank Cantor dynamical systems}
\label{inf_rank_section}

In this section, we give sufficient conditions for a Bratteli diagram of infinite rank to have a prescribed (finite or infinite) number of probability ergodic invariant measures. We define a class of Bratteli diagrams of infinite rank that, in some sense, generalizes the class of Bratteli diagrams of finite rank. A diagram of this class has a prescribed number of uniquely ergodic subdiagrams such that the extension of the unique invariant measure from each subdiagram to the whole diagram is finite. Moreover, there are no other finite ergodic invariant measures for the Bratteli diagram.

\subsection{A class of Bratteli diagrams of infinite rank}\label{inf_rank_main_subsect}
Let us assume that (after telescoping) every level $V_n$ of a Bratteli diagram
$B$ admits a  partition
$$
V_n = \bigcup_{i = 0}^{l_n}V_{n,i}, \ \ n = 1,2,\ldots,
$$
into disjoint subsets $V_{n,i}$ such that $V_{n,i} \neq \emptyset$, for
$i = 1,\ldots,l_{n}$, and $l_n \geq 1$. Moreover, let
$$
L_{n+1} = \{1,\ldots,l_{n+1}\} = \bigcup_{i = 1}^{l_n}L_{n+1}^{(i)},
$$
where $L_{n+1}^{(i)} \neq \emptyset$ and  $L_{n+1}^{(i)} \cap
 L_{n+1}^{(j)} = \emptyset$ for $i \neq j$, $i,j = 1,\ldots, l_n$. Hence,
  for every $j = 1, \ldots, l_{n+1}$, there exists a unique $i = i(j) \in
  \{1,\ldots, l_n\}$ such that $j \in L_{n+1}^{(i)}$. Denote
$$
V_{n+1}^{(i)} = \bigcup_{j \in L_{n+1}^{(i)}} V_{n+1,j}
$$
for $1 \leq i \leq l_n$.

We can interpret the sets $L_n^{(i)}$, defined above, in terms of
 subdiagrams. For this, select a sequence $\ov i = (i_1,i_2, \ldots)$ such
  that  $i_1 \in L_1$, $i_2 \in L_2^{(i_1)}$, $i_3 \in L_3^{(i_2)}, \ldots$
  and  define a subdiagram $B_{\ov i} = (\ov V, \ov E)$, where
  $$
  \ov V = \bigcup_{n=1}^{\infty}V_{n,i_n} \cup \{v_0\}.
$$

Now we formulate conditions $(c1)$, $(d1)$, $(e1)$ which are analogues
of conditions $(c)$, $(d)$, $(e)$ used in  Theorem~\ref{main1}:

\medskip
$(c1)$
$$
\sum_{n = 1}^{\infty} \left(\max_{i \in L_{n}} \max_{v \in V_{n+1}^{(i)}}
 \sum_{w \notin V_{n,i}} f_{vw}^{(n)}\right) < \infty;
$$

\medskip
$(d1)$
$$
\max_{j \in L_{n+1}}\max_{v,v' \in V_{n+1,j}} \sum_{w \in V_{n}}
\left|f_{vw}^{(n)} - f_{v'w}^{(n)}\right| \rightarrow 0 \mbox{ as } n \rightarrow \infty;
$$

\medskip
$(e1)$ for every $v \in V_{n+1,0}$, \\
$(e1.1)$
$$
\sum_{w \in V_n \setminus V_{n,0}} f_{vw}^{(n)} \rightarrow 1
\ \ \mbox{ as } n \rightarrow \infty;
$$

\noindent $(e1.2)$ there exists $C>0$ such that $F_{vi}^{(n)} \leq 1 - C$
for every $i = 1,\ldots,l$, where
$$
F_{vi}^{(n)} = \sum_{w \in V_{n,i}}f_{vw}^{(n)}.
$$
Notice that the new condition (e1.1) in the case of infinite rank Bratteli diagrams is stronger than the corresponding condition (e1) in Theorem~\ref{main1}. 

\medskip

 Let
 $$\mathcal{L} = \{\ov i
  = (i_1, i_2, \ldots) \ : \ i_1 \in L_1, i_{n+1} \in L_{n+1}^{(i_n)}, n = 1,2,
  \ldots\}.
  $$
  We call such a  sequence $\ov i \in   \mathcal{L}$  a \textit{chain}.
 We remark that a Bratteli diagram $B = (V,E)$ of finite rank has the form
 described in this section. We also notice that the following theorem does not require the stochastic incidence matrices of a Bratteli diagram to have the property of regularly vanishing blocks.

\begin{theorem}[\cite{BezuglyiKarpelKwiatkowski2018}]\label{main_gen_case}
Let $B = (V,E)$ be a Bratteli diagram satisfying the conditions $(c1),(d1),
(e1)$. Then:

(1) for each $\ov i \in \mathcal{L}$, any measure $\mu_{\ov i}$ defined on
$B_{\ov i}$ has a finite extension $\wh{\mu}_{\ov i}$ on $B$,

(2) each subdiagram $B_{\ov i}$, $\ov i \in \mathcal{L}$, is uniquely
ergodic,

(3) after normalization, the measures $\wh{\mu}_{\ov i}$, $\ov i \in
\mathcal{L}$, form the set of all probability ergodic invariant measures on
$B$, in particular, $|\mathcal{E}_1(B)| = |\mathcal L|$.
\end{theorem}

\subsection{Examples}

\begin{example}[Pascal-Bratteli diagram]

In Subsection~\ref{inf_rank_main_subsect}, we defined a class of Bratteli diagrams $B =
(V,E)$ such that the set of all ergodic invariant probability
measures coincides with the set $\mathcal{L}$ of all infinite chains $\ov i$. Each
ergodic probability invariant measure $\wh{\mu}_{\ov i}$ is an extension of
a unique invariant measure $\mu_{\ov i}$ from the subdiagram $B_{\ov i}$,
and the sets $X_{B_{\ov i}}$ are pairwise disjoint. It turns out, that the set of ergodic invariant measures for Pascal-Bratteli diagram has a different structure.

For the Pascal-Bratteli diagram, we have $V_n = \{0,1,\ldots, n\}$ for $n = 0,1,
\ldots$, and the entries $\tl f_{ki}^{(n)}$ of the incidence matrix
$\tl F_n$  are of the form
$$
\tl f_{ki}^{(n)} =
\left\{
\begin{aligned}
1 &, \mbox{ if } i = k \mbox{ for } 0 \leq k < n+1,\\
1 &, \mbox{ if } i = k - 1 \mbox{ for } 0 < k \leq n+1,\\
0 &, \mbox{ otherwise}.\\
\end{aligned}
\right.
$$
where $k = 0, \ldots, n+1$,  $i = 0,\ldots,n$ (see~\cite{MelaPetersen2005,Vershik_2011, Vershik_2014,
FrickPetersenShields2017}).
Moreover,
$$
h_i^{(n)} = {n \choose i},
$$
for $i = 0,\ldots,n$. The entries of the corresponding stochastic matrices are
$F_n$:
\begin{equation}\label{Pascal_stochastic}
f_{ki}^{(n)} =
\left\{
\begin{aligned}
\frac{k}{n+1} &, \mbox{ if } i = k -1 \mbox{ and } 0 < k \leq n+1,\\
1 - \frac{k}{n+1} &, \mbox{ if } i = k \mbox{ and } 0 \leq k < n+1,\\
0 &, \mbox{ otherwise},\\
\end{aligned}
\right.
\end{equation}
It is known~(see e.g. \cite{MelaPetersen2005}) that each ergodic invariant
 probability measure has the form $\mu_p$, $0 < p < 1$, where
$$
\mu_p\left(X_i^{(n)}\right) = {n \choose i}p^i (1-p)^{n-i}, \quad i = 0,
\ldots,n.
$$

\begin{proposition}[\cite{BezuglyiKarpelKwiatkowski2018}] For the Pascal-Bratteli diagram,
the set $\mathcal{L}$ of all infinite chains $\overline{i}$ is empty.
\end{proposition}

\end{example}

\begin{example}[A class of Bratteli diagrams with countably many ergodic
invariant measures]\label{subsection_cntbl_erg_meas}
 In this example, we present a class of
 Bratteli diagrams with countably infinite set of ergodic invariant measures. Let $V_n = \{0, 1, \ldots, n\}$ for
$n = 0, 1, \ldots$, and let
$\{a_n\}_{n = 0}^{\infty}$ be a sequence of natural numbers such that
\begin{equation}\label{property1}
\sum_{n = 0}^{\infty} \frac{n}{a_n + n} < \infty.
\end{equation}
  Consider the Bratteli diagram $B$ with $(n+2) \times (n+1)$ incidence matrices
  $$
\tl F_n =
\begin{pmatrix}
a_n & 1 & 1 & \ldots & 1 & 1\\
1 & a_n & 1 & \ldots & 1 & 1\\
\vdots & \vdots & \vdots & \ddots & \vdots & \vdots \\
1 & 1 & 1 & \ldots & 1 & a_n\\
1 & 1 & 1 & \ldots & 1 & a_n\\
\end{pmatrix}.
$$

  Then
 $$
 |r^{-1}(v)| = a_n + n
 $$
for every $v \in V_{n+1}$ and every $n = 0,1, \ldots$

  The Bratteli diagram defined above admits
 an order generating the Bratteli-Vershik homeomorphism (see
  \cite{HermanPutnamSkau1992}, \cite{GiordanoPutnamSkau1995}, or
\cite{BezuglyiKwiatkowskiYassawi2014}, \cite{BezuglyiKarpel2016}). 
 In particular, we can use the so called \textit{consecutive} ordering such that $X_B$ has the unique  minimal infinite path passing through the vertices $0 \in V_n$, $n \geq 0$  and the unique maximal infinite path passing through the vertices $n \in V_n
 $, $n \geq 0$. A Vershik map $\varphi_B \colon X_B \rightarrow
 X_B$ exists and it is minimal. Figure 2 below shows an example of such a 
 Bratteli  diagram. It is known that all minimal Bratteli-Vershik systems with a 
 consecutive ordering have entropy zero (see e.g. \cite{Durand2010}) hence 
 the system that we describe in this subsection has zero entropy.
 
\begin{figure}[ht]
\unitlength = 0.4cm
\begin{center}
\begin{graph}(28,14)
\graphnodesize{0.4}
\roundnode{V0}(1,13.5)
\roundnode{V11}(1,9)
\roundnode{V12}(10,9)

\bow{V11}{V0}{0.18}
\bow{V11}{V0}{-0.21}
\bow{V12}{V0}{0.07}
\bow{V12}{V0}{-0.06}
\freetext(1.1,11){$\cdots$}
\freetext(5.6,11){$\cdots$}

\roundnode{V21}(1,4.5)
\roundnode{V22}(10,4.5)
\roundnode{V23}(19,4.5)

\bow{V21}{V11}{0.18}
\bow{V21}{V11}{-0.21}
\edge{V22}{V11}
\edge{V21}{V12}
\bow{V22}{V12}{0.18}
\bow{V22}{V12}{-0.21}
\edge{V23}{V11}
\bow{V23}{V12}{0.06}
\bow{V23}{V12}{-0.07}
\freetext(1.1,6){$\cdots$}
\freetext(10.1,6){$\cdots$}
\freetext(16,6){$\cdots$}

\roundnode{V31}(1,0)
\roundnode{V32}(10,0)
\roundnode{V33}(19,0)
\roundnode{V34}(28,0)

\bow{V31}{V21}{0.18}
\bow{V31}{V21}{-0.21}
\freetext(1.1,1.5){$\cdots$}
\edge{V31}{V22}
\edge{V31}{V23}

\edge{V32}{V21}
\bow{V32}{V22}{0.18}
\freetext(10.1,1.5){$\cdots$}
\bow{V32}{V22}{-0.21}
\edge{V32}{V23}

\edge{V33}{V21}
\bow{V33}{V23}{0.18}
\freetext(19.1,1.5){$\cdots$}
\bow{V33}{V23}{-0.21}
\edge{V33}{V22}

\edge{V34}{V21}
\bow{V34}{V23}{0.06}
\freetext(24.9,1.5){$\cdots$}
\bow{V34}{V23}{-0.06}
\edge{V34}{V22}

\freetext(1,-1){$\vdots$}
\freetext(10,-1){$\vdots$}
\freetext(19,-1){$\vdots$}
\freetext(28,-1){$\vdots$}

\freetext(10,-2.5){Figure 2}

\end{graph}
\end{center}

\vspace{0.8 cm}
\end{figure}

Denote by $B_i = (W^{(i)}, E^{(i)})$, $i = 0,1,\ldots, \infty$, the
 subdiagrams of $B$ determined by the following sequences of vertices
 (taken consecutively from $V_0$, $V_1$, $\ldots$): for $B_0$,
 $W^{(0)} = (0,0,0, \ldots)$; for $B_i$,   $W^{(i)} = (0,1,\ldots,i-1,i,i,i
 \ldots)$ for $i = 1,2,\ldots$,  and for $B_{\infty}$,  $W^{(\infty)} =(0,1,2,
  \ldots)$. Then each $B_i$ is an odometer and $E^{(i)}$ is the set of all
  edges from $B$ that belong to $B_i$. 
Let $\mu_i$ be the unique invariant (hence ergodic) probability measure on
 the odometer $B_i$. Then each measure $\mu_i$
  can
 be extended to a finite invariant measure $\wh{\mu}_i$ on the diagram $B$
 and it is supported by the set $\wh{X}_{B_i}$ (see \cite{BezuglyiKarpelKwiatkowski2018}). We use the same symbol
 $ \wh{\mu}_i$ to denote the normalized (probability)
measure obtained from the extension of $\mu_i$ for $i = 0,1,\ldots,\infty$.

\begin{proposition}[\cite{BezuglyiKarpelKwiatkowski2018}]\label{countbl_erg_meas}
The measures $\wh{\mu}_i$, $i = 0, 1, \ldots, \infty$, form a set of all
ergodic probability invariant measures on the Bratteli diagram $B = (V,E)$
defined above.
\end{proposition}

\end{example}

\medbreak
{\bf Acknowledgement.} The authors are very thankful to our colleagues and
collaborators for numerous fruitful and stimulating discussions. 
Especially, we would like to 
thank T. Downarowicz, P. Jorgensen, D. Kwietniak, J. Kwiatkowski,  P.~Muhly, P. Oprocha.

\bibliographystyle{alpha}
\bibliography{referencesBKK5}

\end{document}